\documentclass[10pt,a4paper,fleqn]{article}
\usepackage{amsmath}
\usepackage{amssymb}
\usepackage{amscd}
\usepackage{plain}
\usepackage{graphicx}
\usepackage{epsfig}
\usepackage{t1enc,epsfig}
\usepackage[mathscr]{eucal}
\usepackage{epstopdf}
\usepackage{amsthm}
\usepackage{fullpage}

\newtheorem{theorem}{Theorem}

\newtheorem{corollary}[theorem]{Corollary}
\newtheorem{definition}[theorem]{Definition}
\newtheorem{example}[theorem]{Example}
\newtheorem{lemma}[theorem]{Lemma}
\newtheorem{proposition}[theorem]{Proposition}
\newtheorem{remark}[theorem]{Remark}
\newtheorem{conjecture}[theorem]{Conjecture}

\newcommand{\q}{\mbox{q}\,}

\newcommand{\U}{\underset}

\newcommand{\bthm}{\begin{theorem}}
\newcommand{\ethm}{\end{theorem}}
\newcommand{\bd}{\begin{definition}}
\newcommand{\ed}{\end{definition}}
\newcommand{\bs}{\begin{proposition}}
\newcommand{\es}{\end{proposition}}
\newcommand{\bp}{\begin{proof}}
\newcommand{\ep}{\end{proof}}
\newcommand{\be}{\begin{equation}}
\newcommand{\ee}{\end{equation}}

\newcommand{\br}{\begin{remark}}
\newcommand{\er}{\end{remark}}
\newcommand{\bex}{\begin{example}}
\newcommand{\eex}{\end{example}}
\newcommand{\bc}{\begin{corollary}}
\newcommand{\ec}{\end{corollary}}
\newcommand{\bl}{\begin{lemma}}
\newcommand{\el}{\end{lemma}}
\newcommand{\bj}{\begin{conjecture}}
\newcommand{\ej}{\end{conjecture}}

\begin{document}

\def\bx{{\mathbf x}} \def\bi{{\mathbf i}} \def\bk{{\mathbf k}}

\def\uC{{\underline C}} \def\uD{{\underline D}} \def\ug{{\underline g}}
\def\ux{{\underline x}} \def\ui{{\underline i}} \def\uk{{\underline k}}
\def\U0{{\underline 0}} \def\u1{{\underline 1}}

\def\tP{{\tt P}} \def\tp{{\tt p}} \def\tQ{{\tt Q}} \def\tq{{\tt q}}
\def\tR{{\tt R}} \def\tr{{\tt r}}  \def\bq{{\mathbf q}}
\def\bs{{\mathbf s}}

\def\fB{\mathfrak B} \def\fF{\mathfrak F} \def\fM{\mathfrak M}
\def\fV{\mathfrak V} \def\fW{\mathfrak W}\def\fX{\mathfrak X}
\def\bu{\mathbf u}\def\bv{\mathbf v}\def\bx{\mathbf x}\def\by{\mathbf y}

\def\bbP{\mathbb P} \def\hw{h^{\rm w}} \def\hwi{{h^{\rm w}}}
\def\beq{\begin{eqnarray}} \def\eeq{\end{eqnarray}}
\def\beqq{\begin{eqnarray*}} \def\eeqq{\end{eqnarray*}}
\def\rd{{\rm d}} \def\re{{\rm e}} \def\rO{{\rm O}} \def\rw{{\rm w}}
\def\BX{\mathbf{X}}\def\Lam{\Lambda}\def\BY{\mathbf{Y}}

\def\mwe{{D^{\rm w}_\phi}}
\def\bD{{\mathbf D}} \def\ubD{{\underline\bD}}
\def\bbE{\mathbb{E}}
\def\1{{\mathbf 1}} \def\fB{{\mathfrak B}}  \def\fM{{\mathfrak M}}
\def\diy{\displaystyle} \def\bbE{{\mathbb E}} \def\bu{\mathbf u}
\def\BC{{\mathbf C}} \def\lam{\lambda}
\def\bbB{{\mathbb B}} \def\bbD{{\mathbb D}}
\def\bbI{{\mathbb I}} \def\bbJ{{\mathbb J}}
\def\bbR{{\mathbb R}}\def\bbS{{\mathbb S}}

\def\bbZ{{\mathbb Z}}

\def\B1{\mathbf 1}
\def\gam{{\gamma}} \def\om{\omega} \def\Om{\Omega}
\def\eps{{\epsilon}} \def\tht{{\theta}}
\def\veps{\varepsilon} \def\ueps{{\underline\eps}}
\def\uveps{{\underline\veps}} \def\oveps{{\overline\veps}}
\def\umu{{\underline\mu}} \def\omu{{\overline\mu}}
\def\bmu{{\mbox{\boldmath${\mu}$}}}
\def\beps{{\mbox{\boldmath${\veps}$}}}
\def\uveps{{\underline\veps}}

\def\bPhi{{\mbox{\boldmath${\Phi}$}}}

\def\vphi{\varphi}

\def\beal{\begin{array}{l}}
\def\beac{\begin{array}{c}}
\def\beacl{\begin{array}{cl}}
\def\ena{\end{array}}

\def\bcs{\begin{cases}}
\def\ecs{\end{cases}}

\def\bb{{\mathbf b}}
\def\bpp{{\mathbf p}}
\def\alph{\alpha}  \def\bet{\beta}
\def\rA{{\rm A}} \def\rL{{\rm L}} \def\vphi{{\varphi}} \def\bear{\begin{array}{r}}
\def\tf{{\tt f}} \def\tg{{\tt g}}

\def\bx{\bar x}

 \def\cB{\mathcal B}\def\cM{\mathcal M}\def\cC{\mathcal C}
\def\cF{{\mathcal F}}
\def\cJ{{\mathcal J}} \def\cK{{\mathcal K}}
\def\cT{{\mathcal T}} \def\cX{\mathcal X}
\def\bcC{{\mbox{\boldmath${\cC}$}}}

\def\sC{\mathscr C}

\def\ov{\overline}

\def\wh{\widehat} \def\ov{\overline}
\def\wt{\widetilde}
\def\wtk{{\widetilde k}}

\def\cl{\centerline}

\title{Weighted entropy and optimal portfolios\\ for risk-averse Kelly investments}

\author{M.~Kelbert$^{1}$, I.~ Stuhl$^{2-4}$, Y.~Suhov$^{5-7}$  }

\date{}
\maketitle

\begin{abstract} Following a series of works on capital growth investment, we analyse log-optimal portfolios
where the return evaluation includes `weights' of different outcomes. The results are twofold: (A) under certain
conditions, the logarithmic growth rate leads to a supermartingale, and (B) the optimal (martingale) investment
strategy is a proportional betting. We focus on properties of the optimal portfolios and discuss a number of
simple examples extending the well-known Kelly betting scheme.

An important restriction is that the investment does not exceed the current capital value and allows
the trader to cover the worst possible losses.

The paper deals with a class of discrete-time models. A continuous-time extension is a topic of an
ongoing study.
 \end{abstract}

\footnotetext{$^{1}$ Moscow Higher School of Economics, RF}

\footnotetext{$^{2}$ Math Dept, University of Denver, CO, USA; $^{3}$ IMS, Univesity of Sao Paulo,
SP, Brazil; $^{4}$ University of Debrecen, Hungary}

\footnotetext{$^{5}$ Math Dept, Penn State University, PA, USA;  $^{6}$ DPMMS, University of Cambridge,
 UK;  $^{7}$ IPIT RAS, Moscow, RF}

\footnotetext{2010 {\em Mathematics Subject Classification: 60A10, 60B05, 60C05 (Primary), 91G80,
91G99 (Secondary)}}
\footnotetext{{\em Key words and phrases:} weight function, return function, predictable strategy,
expected weighted interest rate, supermartingale, martingale, log-optimal investment portfolio  \par}

\section{A Markovian model with a single risky asset}

This paper is an initial part of a work on log-optimal portfolios influenced by a number of earlier
publications, mainly by T. Cover and co-authors. Cf.
Refs \cite{AC, Co} and \cite{CT}, Chapter 6. Also see \cite{MSZZ, MZZ, Z} and Ref \cite{MTZ}, Parts
II and III. We also intend to use a recent progress in studying weighted entropies; cf.
\cite{SSYK, SS1, SS2, SYS}. A strong impact on the whole direction of this research was made
by \cite{KrS1, KrS2} where a powerful methodology of a convex analysis has been developed (and elegantly presented)
in a general form, leading
-- among other achievements -- to existence of log-optimal portfolios. See Theorem 1 from \cite{KrS2}.
In the present article, we attempt to go beyond the issue of existence and provide a specific form
of the optimal strategy.

Let us discuss a finance-related context of this work.
The sequential version of portfolio selection problem has received much attention
in the literature not to speak about the financial practice, see \cite{MSZZ, MZZ, MTZ,Z} and the references therein.
A simple discrete-time model of a wide use in financial
engineering is where the market consists of one riskless asset and one or more risky assets.
(If the riskless asset produces a zero return, we can speak of risky assets only.)
Investments are made at times $n-1=0,1,\ldots$; the returns are recorded at subsequent
times $n=1,2,\ldots$.

We consider two investment schemes, showing that the results are valid for both schemes
{\it mutatis mutandis}. \\

Scheme I: an investor signs a deal with a broker at the time $n-1$ but the actual
transaction happens at the moment $n$ when the betting results become available.\\

Scheme II: at the moment $n-1$ an investor transfers the required capital to a
broker who invests this capital to buy shares or other risky assets.\\

In fact, Scheme II  can be treated as a version of Scheme I, where the number of assets
increases by 1.

For convenience of presentation we switch freely between these two schemes
keeping in mind that the results are valid in both cases with minor changes.

The mathematical problem under consideration is to characterize an optimal investment
strategy/po\-li\-cy/portfolio (it will be convenient to use all these terms as equal in right, as
with some other synonyms).
Formally, we have to introduce the objective function and  describe the class of strategies
within which the optimization problem is solved. The setting for our study is probabilistic:
we assume that, generally speaking, the returns are {\it random}, with known probabilistic laws.
(In practice, it means that probabilistic/statistical features can be established, e.g., from available
historic data.) More precisely, we deal with a {\it random process} of return values. For
illustrative purposes we adopt for the most part a {\it Markovian} model of the return process(es) but also
provide a mathematical result under quite general assumptions. (A number of features
of the solution depends on the character and parameters of the return process.)

The objective function is introduced as the expected value of the {\it weighted logarithmic
return} $S_n$ after $n$ trials. This is in line with the proposal going back to J. Kelly (Kelly investments, \cite{Kelly})
although our approach is based on some important modifications. In particular,
we consider only {\it cautious/risk-averse} investment policies with a guaranteed cover of all possible losses.
Here, a passive/restrained $0$-strategy is a notable example, where the investor decides not
to bet on the outcome of next trial. Pictorially speaking, the answer emerges from a comparison
of the best adventurous bet and the zero-bet (and in most cases the $0$-bet is preferable or the comparison
is inconclusive, and the theory does not produce a formal answer).

In short, the results offered here can be summarized as follows: {\it under certain conditions,
an optimal policy is to invest a proportion of the current capital, regardless of the value
of the capital achieved by the time of the investment decision}. The proportion depends on the
current state of the return process (and possibly on its history). As was said, in many cases
the recommended proportion is $0$ meaning no investment into a given trial.

This paper focuses on discrete-time models. A continuous-time version of our approach is
currently under investigation and will be a topic of forthcoming publications.
In this regard, we note an alternative approach propagated in \cite{BPS} and, more recently,
in \cite{Puh}. An extensive bibliography of this field is available in \cite{DL}.

\vskip .5 truecm

We start with a basic discrete-time setting under Scheme I.
An investor is betting on results $\veps_n$ of subsequent random trials, $n=0,1,2,\ldots$.
Suppose that the $\veps_n$ are generated by a Markov chain with a finite or countable
state space $M$. The transition matrix is $\tP=(\tp (i,k),\,i,k\in M)$.

Let us introduce a {\it return function} $(i,k)\in M\times M\mapsto g(i,k)$ with real values. Then let us
agree that if the investor stakes
$\,c$ on the $n$th trial he/she wins $ \$\,c g(i,k)$ if the result is $k$ preceded by the outcome $i$
of the $(n-1)$st trial (which you may know). (If $g(i,k)<0$ you loose $\,c|g(i,k)|$.) When $g(i,k)=g(k)$,
the return from the next trial is determined by the coming outcome regardless of the previous one(s).
A more general setting, with a `long' memory, is considered in the forthcoming sections.
Let $Z_{n-1}$ be the investor's fortune/capital
after $n-1$ trials; we will assume that $Z_{n-1}>0$ which will be justified below. ($Z_0>0$ is
an initial capital which can be dependent on $\veps_0$, the initial-trial result.)  Introduce the
$\sigma$-algebras $\fW_0 =
\sigma (Z_0,\veps_0)$ and $\fW_n =\fW_0\vee \sigma (\beps_0^n)$, $n\geq 1$, and consider a
sequence  $\{C_n,\,n\geq 0\}$ where $C_n$ is $\fW_n$-measurable random variable (a predictable strategy):
$C_n=C_n(Z_0,\beps_0^n)$. Here and below, $\beps_0^n=(\veps_0,\ldots ,\veps_n)$
stands for a string of subsequent random variables (RVs) representing the first $(n+1)$ trial results.
The dependence upon $Z_0$ will be omitted. The recursion for \ $Z_n$ \ is
\beq\label{eq13B1}
Z_n=Z_{n-1}+C_{n-1}g(\veps_{n-1},\veps_n)=Z_{n-1}\left[1+\frac{C_{n-1}g(\veps_{n-1},\veps_n)
}{Z_{n-1}}\right];\eeq
it shows that $Z_n\in\fW_n$: $Z_n=Z_n(\beps_0^n )$. For the return function $g$ we will use the acronym RF.

Next, let us consider another function, $(i,k)\mapsto \vphi (i,k)\geq 0$,
representing a `utility' value assigned to outcome $k$ when it succeeds outcome $i$.
If $\vphi (i,k)\equiv 1$, all outcomes are treated entirely in terms of their returns, and if $\vphi (i,k)$
does not depend on $i$, the value does not take into account the history. We say that $\vphi$
is a {\it weight function} (WF); including one-step history $i$ agrees with the Markovian assumption
for $\veps_n$.

It will be always assumed that the WF and the RF are not identically $0$; this includes the modified models
below with more than one asset (where we deal with an RF vector $\ug$).

We wish to maximize, in $C_0,\ldots ,C_{n-1}$, the mean value $\bbE S_n$ where
\beq\label{eq14B1}
S_n:=\sum\limits_{j=1}^n\vphi (\veps_{j-1},\veps_j)\ln\,\frac{Z_j}{Z_{j-1}},
\eeq
and determine, when possible, a sequence of optimal strategies $\{C_j^\rO\}$,
within `natural' classes $\{\cC_j\}$ of predictable strategies defined by recursive
inequalities \eqref{eqa1-3.3} below:
\beq\label{eq131}(C^\rO_0,\ldots ,C_{n-1}^\rO)={\rm{argmax}}\;\Big[\bbE S_n:\;C_j\in\cC_j,\;0\leq j\leq n-1\Big].\eeq
The classes $\cC_j$ are described through conditions (a0)--(a2) or (a0)--(a3) listed in Eqn \eqref{eqa1-3.3} below.
Under our assumptions, the optimum is at a proportional betting, where
$C^\rO_{j-1}=D^\rO_{j-1}(\veps_{j-1})Z_{j-1}$. Here $Z_{j-1}$ is the capital after $j-1$ trials
and $D^\rO_{j-1}(i)$ is the proportionality coefficient indicating the fraction of the capital to be invested
into the $j$th trial.

Quantity $S_n/n$ can be considered as a weighted log-capital rate after $n$ trials. When
$\vphi (i,k)\equiv 1$, the sum in \eqref{eq14B1} becomes telescopic, and $S_n$ equals \
$\ln\diy\frac{Z_n}{Z_0}$ \ (a standard quantity in financial calculations, particularly in relation
to the Kelly-type investments).

The form of summation in Eqn \eqref{eq14B1} suggests the use of the weighted Kullback--Leibler (KL)
entropy of the row probability vector $\big(\tp (i,k),\,k\in M\big)$ relative to chosen `calibrating' functions
$(i,k)\in M\times M\mapsto \tq_j(i,k)>0$, $j=0,1,\ldots $. To this end, set:
\beq\label{eq17B1}
\alpha_j(i)=\sum\limits_{l\in M}\vphi (i,l)\tp (i,l)\ln\,\frac{\tp (i,l)}{\tq_j(i,l)},\;i\in M,\;j\geq 0 .\eeq
For the definition and basic properties of weighted entropies, see \cite{SSYK}
and references therein. Some applications of weighted entropies are discussed
in \cite{SS1, SS2, SYS}.

The sum in \eqref{eq17B1} and similar series below are supposed to converge absolutely.

The choice of calibrating functions (CFs) $\tq_j (i,k)$ is a part of the optimization procedure
and is discussed below: see \eqref{eqa1-3.3} and \eqref{eqq11Rl}. We consider
the random process (RP) of the cumulative weighted KL entropy
\beq\label{eq101M}\rA_n=\sum\limits_{j=1}^n\alpha_{j-1}(\veps_{j-1}),\;\hbox{
 with }\;\bbE\rA_n=\sum\limits_{j=1}^n\bbE\alph_{j-1} (\veps_{j-1}).\eeq
Also fix a value $b>0$ (a proportional ruin threshold).

Let us summarize conditions on the class of policies and involved functions:
\ $\forall$ \ $j\geq 0$,
\beq\label{eqa1-3.3}\beal{\rm{(a0)}}\;C_j\in\fW_j,\;\; {\rm{(a1)}}\;0\leq C_j<
Z_j,\;{\rm{(a2)}}\;1+\diy\frac{C_jg(\veps_j,\veps_{j+1})}{Z_j}\geq b,\;\hbox{ and}\\
{\rm{(a3)}}\quad C_j(\beps_0^j)\sum\limits_{l\in M}\vphi (\veps_j,l)\tq_{j+1}(\veps_j,l)
g(\veps_j,l)=0.\ena\eeq
% \;\hbox{ which can be related to}\\
% {\rm{(a4)}}\;\sum\limits_{l\in M}\vphi (i,l)\tq_n(i,l)g(i,l)=0,\;\hbox{ and }\;{\rm{(a5)}}\;C_{n-1}
% (\bi_0^{n-1})=0,\;\hbox{ for some or all }\\
% \qquad i\in M,\;\bi_0^{n-1}=(i_0,\ldots ,i_{n-1})\in M^n.\ena\eeq
Also, \ $\forall$ \ $i\in M,\;j\geq 0$, we assume the condition labelled as (q--p) in
Eqn \eqref{eqq11Rl}:
\beq\label{eqq11Rl} {\rm{(q-p)}}\qquad
\sum\limits_{l\in M} \vphi (i,l)\big[\tq_j (i,l)-\tp (i,l)\big]\leq 0.\eeq

\vskip .5 truecm

{\bf Remark 1.1.} Recall, $C_j=0$ means no investment into the result of the $(j+1)$st trial.
Note that $C_j\equiv 0$ is always a feasible choice: it yields $S_{j+1}=S_{j}$. In some situations
it gives an optimum (when the outlook of the results is not favorable for the investor).
On the other hand, if we manage to verify that, \ $\forall$ \ $j$, the sum
$\sum\limits_{l\in M}\vphi (\veps_j,l)\tq_{j+1}(\veps_j,l)g(\veps_j,l)=0$ then
(a3) provides no additional restriction upon $C_j$ and can be discarded. It means that
the optimality can be achieved within the larger class of strategies satisfying (a0)--(a2) only.
(Still, the optimum may be $C_j=0$.)
This aspect of the theory will be repeatedly stressed in various models and examples below.

\vskip .5 truecm

We offer the following result.
\vskip .5 truecm

{\bf Theorem 1.1.} {\sl Suppose the recursion \eqref{eq13B1} holds true.

{\rm{(a)}} Suppose a sequence of CFs $\tq_n$ is given, obeying
\eqref{eqq11Rl}. Take any sequence $\{C_n,\;n\geq0\}$ of random variables (RVs) $C_n$
satisfying properties {\rm{(a0)}}--{\rm{(a3)}} in Eqn \eqref{eqa1-3.3}. Consider RVs $S_n$
and $A_n$ defined in \eqref{eq14B1} and \eqref{eq17B1}--\eqref{eq101M}.
Then the sequence of differences $\{S_n-\rA_n,\,n\geq 1\}\;$ is a supermartingale; consequently,
$\bbE\,S_n\leq\bbE\rA_n$ \ $\forall$ $n\geq 1$. \vskip .3 truecm

{\rm{(b)}} To achieve equality $\bbE\,S_n=\bbE \rA_n$: the sequence $\{S_n-\rA_n\}$
is a martingale for a sequence of RVs $C_n$ satisfying {\rm{(a0)}}--{\rm{(a3)}} in
\eqref{eqa1-3.3} iff the following conditions {\rm{(i)}}--{\rm{(ii)}} hold.

\quad {\rm{(i)}} There exists a function $D:\, M\to\bbR$ such that, $\forall$ \ $i,k\in M$,
\beq\label{eq8M1}\beal \quad {\rm{(i1)}}\;\;0\leq D(i)< 1,\;\; {\rm{(i2)}}\;\; 1+D(i)g(i,k)\geq b,\;\;
{\rm{(i3)}}\;\; D(i)\sum\limits_{l\in M}
{\diy\frac{\tp (i,l)\vphi (i,l)g(i,l)}{1+D(i)g(i,l)}}=0,\\
\qquad\quad\hbox{i.e., either}\;{\rm{(i3A)}}\;\sum\limits_{l\in M}
{\diy\frac{\tp (i,l)\vphi (i,l)g(i,l)}{1+D(i)g(i,l)}}=0\;\hbox{ or }\;{\rm{(i3B)}}\;D(i)=0,\;\hbox{ and}\\
\quad{\rm{(i4)}}\;\;\hbox{the CFs \ $\tq_j$ \ are of the form }\;\tq_j(i,k)=
\diy\frac{\tp (i,k)}{1+D(i)g(i,k)},\;\hbox{independently of}\\
\qquad\quad j=0,1,\ldots .\ena\eeq

\quad {\rm{(ii)}} $\forall$ \ $n\geq 1$, \ the policy $C_{n-1}$ produces a proportional
investment: $C_{n-1}(\beps_0^{n-1})=D(\veps_{n-1})Z_{n-1}$. \medskip

Furthermore, the CF values
$\tq_n(i,k)=\tq (i,k)$ given in {\rm{(i4)}} satisfy $\sum\limits_{l\in M}\tq (i,l)=1$ (which yields
transition probability matrices) iff, in addition to {\rm{(i3)}}, we have that
\beq\label{eqMBa}D(i)=0\;\;\hbox{ or }\;\;\sum\limits_{l\in M}\frac{\tp (i,l)g(i,l)}{1+D(i)g(i,l)}=0,
\;\;i\in M.\eeq

{\rm{(c)}} Define the map $i\in M\mapsto D^\rO(i)$ as follows. Given $i$, consider Eqn {\rm{(i3A)}}:
it has at most one solution $D(i)>0$. If {\rm{(i3A)}} has a
solution $D(i)>0$ obeying conditions {\rm{(i1)}}--{\rm{(i2)}}, set $D^\rO (i)=D(i)$; otherwise
$D^\rO (i)=0$. Then the policy $C^\rO_{n-1}=D^\rO(\veps_{n-1})Z_{n-1}$ yields
the following value $E_n$ for the expectation $\bbE S_n$:
\beq\label{eqOptRate1}E_n=\sum\limits_{j=1}^n\bet_{j-1}\;\hbox{ where }\; \bet_{j-1}
=\bbE\Big\{\vphi (\veps_{j-1},\veps_j)
\ln\,\Big[1+D^\rO(\veps_{j-1})g(\veps_{j-1},\veps_j)\Big]\Big\}.\eeq
Moreover, the value $E_n$ gives the maximum of \ $\bbE S_n$ \ over
all strategies satisfying conditions {\rm{(a0)}}--{\rm{(a3)}} in \eqref{eqa1-3.3}.

{\rm{(d)}} Suppose that the map $i\in M\mapsto D^\rO(i)$ from assertion {\rm{(c)}} is such that
$D^\rO(i)>0$ (so the alternative {\rm{(i3A)}} holds), $\forall$ $i\in M$.
Then the policy $C^\rO_{n-1}=D^\rO(\veps_{n-1})Z_{n-1}$ maximises
each summand $\bet_{j-1}$ in \eqref {eqOptRate1}, and therefore yields the maximum
of the whole sum \ $\bbE S_n$, among strategies satisfying properties {\rm{(a0)}}--{\rm{(a2)}}
in Eqn \eqref{eqa1-3.3}.}

\vskip .5 truecm

{\bf Proof.}  (a) Write:
$$\beal\bbE\Big\{\big(S_n-\rA_n\big)\big|\fW_{n-1}\Big\}=S_{n-1}-\rA_{n-1}\\
\qquad +\diy\bbE\left\{\left[
\vphi (\veps_{n-1},\veps_n)\ln\,\left(1+\diy\frac{C_{n-1}g(\veps_{n-1},\veps_n)}{Z_{n-1}}
\right)\right]\Big|\fW_{n-1}\right\}-\alph_{n-1}(\veps_{n-1}).\ena$$
Next, represent
\beq\label{eq9*}
\beal\diy\bbE\left\{\left[\vphi (\veps_{n-1},\veps_n)\ln\,\left(1+\frac{C_{n-1}
g(\veps_{n-1},\veps_n)}{Z_{n-1}}\right)\right]\Big|\fW_{n-1}\right\}-\alph_{n-1} \\
\quad\diy =\;\sum\limits_{l\in M}\vphi_n(\veps_{n-1},l)\;\tp (\veps_{n-1},l)\;\ln\;\left[1
+\frac{C_{n-1}g(\veps_{n-1},l)}{Z_{n-1}}\right]\\
\qquad\qquad\diy - \sum\limits_{l\in M}\vphi (\veps_{n-1},l)\tp (\veps_{n-1},l)\ln\,
\frac{\tp (\veps_{n-1},l)}{\tq_n(\eps_{n-1},l)}\\
\quad\diy = \sum\limits_{l\in M}\vphi(\veps_{n-1},l)\tp (\veps_{n-1},l) \ln\frac{1
+C_{n-1}g(\veps_{n-1},l)/Z_{n-1}}{\tp (\veps_{n-1},l)/\tq_n(\veps_{n-1},l )}\\
\quad\diy = \sum\limits_{l\in M}\vphi (\veps_{n-1},l)\;\tp (\veps_{n-1},l)\;
\ln\;\frac{h_n(\veps_{n-1},l)}{\tp (\veps_{n-1},l)}\\
\quad\diy \leq\sum\limits_{l\in M}\vphi (\veps_{n-1},l)\tp (\veps_{n-1},l)
\left[ \frac{h_n(\veps_{n-1},l)}{\tp (\veps_{n-1},l)} - 1 \right]{\mathbf 1}\big(\tp
(\veps_{n-1},l)>0\big)\\
\quad\diy =\sum\limits_{l\in M}\vphi (\veps_{n-1},l) \Big[h_n(\veps_{n-1},l)
- \tp (\veps_{n-1},l)\Big]\leq 0.\ena\eeq
Here $h_n(\veps_{n-1},k) := \diy\tq_n(\veps_{n-1},k )\left[1 +\frac{C_{n-1}g(\veps_{n-1},k )
}{Z_{n-1}}\right]$, $k\in M$.

The final inequality in \eqref{eq9*} holds since
\beq\label{eq9.1*}\beal
\diy \sum\limits_{l\in M}\vphi (\veps_{n-1},l)h_n(\veps_{n-1},l ) = \sum\limits_{l\in M}
\vphi (\veps_{n-1},l)\tq_n(\veps_{n-1},l )\\
\qquad+\diy\frac{C_{n-1}}{Z_{n-1}}\sum\limits_{l\in M}\vphi (\veps_{n-1},l)\tq_n(\veps_{n-1},l)
g(\veps_{n-1},l)\leq
\sum\limits_{l\in M}\vphi (\veps_{n-1},l)
\tp (\veps_{n-1},l),\ena\eeq
because of property (a3) in Eqn \eqref{eqa1-3.3} and (q--p) in Eqn \eqref{eqq11Rl}.

(b) To guarantee the martingale property, we have to reach equality in the inequalities in
\eqref {eq9*}. The first inequality becomes equality iff $h_n=\tp$ which yields (i4) in Eqn
\eqref{eq8M1}. The second inequality (i.e., the inequality in \eqref{eq9.1*}) gives equality
iff (I) the relation (a3) from \eqref{eqa1-3.3} holds true (for $C_n(\beps_0^n)=D(\eps_n)Z_n$),
and (II) the bound (q--p) in Eqn \eqref{eqq11Rl} becomes equality. After the substitution from
(i4), both properties (I) and (II) are equivalent to
(i3) from Eqn \eqref{eq8M1}. The inequalities (i1) and (i2) in Eqn \eqref{eq8M1} are the same
as (a1) and (a2) in Eqn \eqref{eqa1-3.3}.

The deduction of Eqn \eqref{eqMBa} is direct, from the above considerations.

(c) The proof of this assertion is straightforward.

(d) Assuming, for a map $i\in M\mapsto D(i)$, the alternative (i3A) in \eqref{eq8M1}
means that, with CFs $\tq_n$ defined by (i4), $\forall$ sequence of policies
$C_n$ satisfying (a0)--(a2) we also have property (a3) and therefore the supermartingale
inequality $\bbE\,S_n\leq\bbE\rA_n$, $n\geq 1$. Setting $C_n=D(\veps_{n-1})Z_n$ yields a martingale equation
$\bbE S_n=E_n$ where $E_n$ is as in \eqref{eqOptRate1}.

Consider the optimization problems emerging from Eqn \eqref{eqOptRate1}:
\beq\label{eq8M101}\beal
{\rm{maximise\;the\;objective\;function}}\quad{\ov\bet}_{j-1}:=\bbE\Big\{\vphi (\veps_{j-1},\veps_j)
\ln\,\Big[1+d(\veps_{j-1})g(\veps_{j-1},\veps_j)\Big]\Big\}\\
\hbox{in variables}\qquad d(i)\in\bbI (i),\;i\in M,\\
\qquad\qquad\hbox{where the feasibility interval }\;
\bbI(i)=\left\{{\wt d}\in [0,1]:\;\;1+{\wt d}g(i,k)\geq b\;\forall\;k\in M\right\}.
\ena\eeq
Here $j=1,\ldots ,n$ (which corresponds with the summand $\bet_{j-1}$ in \eqref{eqOptRate1}).

Observe that the Hessian matrix $H ({\ov\bet}_{j-1})=\Big(\partial^2\bet_{j-1}\big/
\partial d(i)\partial d(i')\Big)$ is non-positive definite:
$$H ({\ov\bet}_{j-1})=-{\rm{diag}}\left(\bbE\left\{\B1 (\veps_{j-1}=i)\frac{\vphi (i,\veps_j)
\big[g(i,\veps_j)\big]^2}{\big[1+d(i)g(i,\veps_j)\big]^2}\right\},\,i\in M\right).$$
It implies that each ${\ov\bet}_{j-1}$ is a concave function in $d(i)\in\bbI (i)$, $i\in M$.
The values $d(i)=D(i)$, $i\in M$, give a zero of the gradient vectors $\nabla{\ov\bet}_{j-1}$
of functions ${\ov\bet}_{j-1}$, hence the maximum for the terms
${\ov\alph}_{j-1}$ from Eqn \eqref{eqOptRate1}. If, in
addition, $D(i)\in\bbI (i)$ \ $\forall$ $i\in M$ (which is what we assume in this statement) then
$\forall$ $n\geq 1$ the sum $\bbE\rA_n$ attains the maximum value among the sequences
$\{C_n\}$ satisfying (a0)--(a2). $\Box$

\vskip .5 truecm

{\bf Remarks. 1.2.} Terminologically, the requirement (q--p) in Eqn \eqref{eqq11Rl}
is referred to as a weighted\\ q,p-dominance condition. As was
stressed earlier, condition (a0) in \eqref{eqa1-3.3} means predictability. The inequalities
(a1) and (a2) are called, respectively, a sustainability condition and a no-ruin condition.
(In fact, the bound (a1) is not used in the proof of statement (a), only (a2).) The relation
(a3) is called the (q,g)-balance condition.
In Eqn \eqref{eq8M1}: (i1) is called a D-sustainability condition, and (i2) a D-no-ruin
condition ((i1) and (i2) emerge from (a2) and (a3), respectively). Next, (i3)
is called a WE (D,g)-balance condition (it is related to (a3)) and (i3A) gives its strong
version. Here WE stands for weight-expected as the condition incorporates
both $\vphi$ and $\tp$. Finally, (i4) is referred to as a $\tq$-representation.
An interpretation of (a1)--(a2)
and (i1)--(i2) is that the investment is controlled by a risk-averse trader
who does not allow the value of the current capital to be below a given level and operates
strictly within the agreed liquidity limits.

{\bf 1.3.} Eqn (i3) in \eqref{eq8M1} (particularly, the version (i3A)) is, arguably, the most
serious restriction from the point of
applications; in all examples under consideration it is a condition that specifies the optimum.
Computationally, it requires solving (a system of) equations involving rational functions.
See examples below.

{\bf 1.4.} The return function $g$ (and, in fact, the weight function $\vphi$) can also be
a part of an optimization procedure, with obvious (although tedious) changes in assertions
(c,d) of Theorem 1.1.

{\bf 1.5.} As can be seen from the proof, a key role in Theorem 1.1(a) is played by
the Gibbs inequality for a weighted KL entropy. In fact, the line of the argument is
essentially independent of particulars of
adopted settings and is repeated throughout the paper.
\vskip .5 truecm

{\bf Example 1.1: IID trials.}  In the case of IID trials, $\tp (i,k)=\tp (k)$
does not depend on $i$. Assume that $g(i,k)=g(k)$,
$\vphi (i,k)=\vphi (k)$ and $\tq_j(i,k)=\tq (k)$. Also choose $b\in (0,1)$. We obtain the following
form of relations
\eqref{eqa1-3.3} and \eqref{eq8M1}:
\beq\label{eqa1-3.4}\beal \ \forall \ j\geq 0:\;\;
{\rm{(a0)}}\;C_j\in\fW_j,\;\; {\rm{(a1)}}\;0\leq C_j<
Z_n,\;{\rm{(a2)}}\;1+\diy\frac{C_j(\beps_0^j)g(\veps_{j+1})}{Z_j(\beps_0^j)}\geq b,\\
\qquad\qquad\quad\,\hbox{and }\;{\rm{(a3)}}\quad\hbox{if the sum}\;\sum\limits_{l\in M}\vphi (l)\tq (l)
g(l)\neq 0\;\hbox{ then }\;C_j(\beps_0^j)=0; \\
{\rm{(q-p)}}\qquad
\sum\limits_{l\in M} \vphi (l)\big[\tq (l)-\tp (l)\big]\leq 0;\ena\eeq
and
\beq\label{eq8B1}\beal\hbox{(i) \ $\exists$ \ (i1) a constant $D\in [0,1]$ such that }\;
{\rm{(i2)}}\;\; 1+Dg(k)\geq b,\;\;k\in M,\\
{\rm{(i3)}}\;\; D\sum\limits_{l\in M}
{\diy\frac{\vphi (l)g(l)\tp (l)}{1+Dg(l)}}=0,\;
\hbox{ and }\;\;{\rm{(i4)}}\;\; \tq (k)=\diy\frac{\tp (k)}{1+Dg(k)},\;
k\in M,\\
\qquad\hbox{with the alternatives (i3A) }\;\;\sum\limits_{l\in M}
{\diy\frac{\vphi (l)g(l)\tp (l)}{1+Dg(l)}}=0\;\hbox{ and (i3B) }\;D=0.\ena\eeq

In the IID case, the weighted KL entropy $\alpha (i)=\alpha$ does not depend on $i$:
$\alph =\sum\limits_{l\in M}\vphi (l)\tp (l)\ln\,\diy\frac{\tp (l)}{\tq (l)}$, see \eqref{eq17B1}.
Consequently, if $\alph\leq 0$ for some CF $\tq$ satisfying (q--p) in \eqref{eqa1-3.4}
then the risk-averse trader would restrain from investments.

Define
\beq\label{eqDOB}D^\rO=\bcs D,&\hbox{if $D\in (0,1)$ is the solution to (i3A) which
satisfies (i2),}\\
0,&\hbox{if $D=0$ is the only solution to (i3) satisfying (i1)--(i2).}\ecs\eeq
According to Theorem 1.1(c), the policies $C^\rO_j=D^\rO Z_j$, $j\geq 0$,
give the maximum of $\bbE S_n$, $n\geq 1$, among the strategies
satisfying (a0)--(a3). Moreover, according to Theorem 1.1(d), if $D^\rO>0$ then
the policies $C^\rO_j=D^\rO Z_j$ maximise $\bbE S_n$ over the strategies satisfying
(a0)--(a2). The maximal value in both cases is $E_n=n\sum\limits_{l\in M}\vphi (l)\tp(l)\ln\,
[1+D^\rO g(l)]\geq 0$.

To illustrate further,
take the case $M=\{0,1\}$ (two outcomes), with $k=0,1$ and probabilities $\tp (0),\tp (1)$.
Here the martingale property occurs iff we can find a constant $D$ such that
% \beq\label{eq8B3}
$$\beac {\rm{(i1)}}\; D\in [0,1],\;{\rm{(i2)}}\; 1+Dg(k)\geq b,\;
{\rm{(i3)}}\;\diy D\left[\frac{\vphi (1)g(1)\tp (1)}{1+Dg(1)}+\frac{\vphi (0)g(0)\tp (0)}{1+Dg(0)}
\right]=0.\ena$$
Next, Eqn (i3) is solved by
$$D=-\frac{\vphi (1)g(1)\tp (1)+\vphi (0)g(0)\tp (0)}{g(0)g(1)[\vphi (1)\tp (1)+\vphi (0)\tp (0)]}
=-\frac{\bbE\big[\vphi (\veps_1)g(\veps_1)\big]}{g(0)g(1)\bbE\vphi (\veps_1)}
\;\hbox{ which should obey (i1), (i2).}$$
To simplify, take $g(1)=-g(0)=\gam >0$. Also, set $\vphi (1)=\vphi (0)=1$ (no preference).
A general form of the WF can be easily incorporated in the argument, which will be also true in
all examples below.
Then the optimal policy is $C^\rO_n =D^\rO Z_n$ where
\beq\label{eqDO2B}D^\rO=\bcs\diy\frac{\tp (1)-\tp (0)}{\gam},&\hbox{if }\;\diy
\frac{b}{2}\leq \tp (0)\leq\tp (1)\hbox{ and }\gam\geq \tp(1)-\tp (0),\\
0,&\hbox{otherwise.}\ecs\eeq

Formula \eqref{eqDO2B} defines a typical Kelly investment scheme \cite{Kelly}.
The maximal growth rate for $\bbE S_n$ takes the form
$$\beal E_n=n\B1\left(\diy
\frac{b}{2}\leq \tp (0)\leq\tp (1)\leq\gam +\tp (0)\right)\Big\{\tp (1)\ln\big[1+\tp (1)-\tp (0)\big]
+\tp (0)\ln\big[1-\tp (1)+\tp (0)\big]\Big\}.\ena$$
If $D^\rO>0$, the optimality holds over the strategies $\{C_j\}$ with properties
(a1,2) $0\leq C_j\leq (1-b)Z_j$, $\forall$ $j\geq 0$. If $D^\rO=0$ then the optimality
holds among the strategies satisfying, in addition, the property (a3) that
$C_j(\beps_0^j)\big[\tp (1)-\tp (0)\big]=0$. (Note that (a3) yields no restriction when,
for instance, $\tp(1)=\tp (0)=1/2$.) More generally, if $g(1)=\gam_1>0, g(0)=-\gam_2<0$
then in the top line of \eqref{eqDO2B} we obtain
$$D^\rO=\frac{p(1)}{\gamma_2}-\frac{p(0)}{\gamma_1},\;\;\hbox{ if \
$\diy\frac{p(1)}{\gam_2}-\frac{p(0)}{\gam_1}\in [0,1)$
and $\diy 1-p(1)+\frac{\gam_2}{\gam_1}p(0)\geq b$,}$$
with obvious changes in the maximal growth rate.
\vskip .5 truecm

{\bf Example 1.2: A two-state Markov chain.} In the Markov case, when the trader observes the
current state $i$, he/she uses the similar optimization procedure for the $i$th row of the
$2\times 2$ transition matrix $\tP=(\tp (i,k))$. Again suppose for simplicity that $M=\{0,1\}$,
the WF $\varphi (i,j)\equiv 1$ and the RF $g$ has $g(1)=-g(0)=\gam >0$. Also suppose
that $b\in (0,1)$ is given. Then an analog of the previous picture emerges. Namely, set,
for $i=1,0$,
\beq\label{eqCondtM}D^\rO(i)=\bcs
\diy\frac{\tp (i,1)-\tp (i,0)}{\gam},&\hbox{if }\;\diy
\frac{b}{2}\leq \tp (i,0)\leq\tp (i,1),\hbox{ and }\gam\geq \tp(i,1)-\tp (i,0),\\
0,&\hbox{otherwise.}\ecs\eeq
The policy $C^\rO_n=D(\eps_n)Z_n(\beps_0^n)$ is optimal, under similar constraints. That is,
if $D^\rO(i)>0$ for both $i=0,1$ then the maximum is attained over strategies satisfying (a0)
$C_j\in\fW_j$ and (a1,2) $0\leq C_j\leq (1-b)Z_j$, \ $\forall$ $j\geq 0$. Otherwise, if $D^\rO(i)=0$
for some $i$ then it is among the strategies obeying (a0)--(a1,2) plus property (a3):
$C_j(\beps_0^j)\big[\tp (\veps_n,1)-\tp (\veps_n,0)\big]=0$ \ $\forall$ $j\geq 0$.

Viz., assume that the Markov chain (MC) is in the stationary regime, with stationary probabilities
$\pi (1)$, $\pi (0)$. Then the maximal growth of $\bbE S_n$ is
$$\beal E_n=n\Bigg(\pi (1)\B1\left(\diy\frac{b}{2}\leq \tp (1,0)\leq\tp (1,1)\leq\gam +\tp (1,0)\right)\\
\qquad\times\Big\{\tp (1,1)\ln\Big[1+\tp (1,1)-\tp (1,0)\Big]
+\tp (1,0)\ln\Big[1-\tp (1,1)+\tp (1,0)\Big]\Big\}\\
\qquad\qquad\qquad +\pi (0)\B1\left(\diy\frac{b}{2}\leq \tp (0,0)\leq\tp (0,1)\leq\gam +\tp (0,0)\right)  \\
\qquad\quad\times\Big\{\tp (0,1)\ln\Big[1+\tp (0,1)-\tp (0,0)\Big]
+\tp (0,0)\ln\Big[1-\tp (0,1)+\tp (0,0)\Big]\Big\}\Bigg).\ena$$

The last observation can be extended to general MCs. Indeed, suppose the trial MC
starts with an invariant distribution $(\pi (i),\,i\in M)$. In this case the statements (c,d)
of Theorem 1.1 assert that the maximum for \ $\bbE S_n$ is given as \
$n\sum\limits_{i,k\in M}\pi_i\tp(i,k)\vphi (i,k)\ln\,\Big[1+D^\rO(i)g(i,k)\Big]$. Note that
we do not need assumptions of irreducibility or aperiodicity: the invariant distribution
is not assumed to be unique.

\section{Markov trials with a riskless asset}

In this section we switch from Scheme I to II. Consider the situation where the trials affect
% the result $\veps_n$ of a trial involves
several assets/returns, say, two, described by RFs
$g^{(0)}$ and $g$. Assume that the return $g^{(0)}$ is riskless:  $g^{(0)}(i,k)=1+\rho$ where $\rho >0$
represents the interest rate, while the other asset, with a RF $g$, is risky. As above, suppose
that the trial results $\veps_0$, $\veps_1$, $\ldots$ \ form an MC with a finite/countable
state space $M$ and a transition probability matrix
$\tP=\big(\tp_{ik}=\tp (i,k),\,i,k\in M\big)$. Fix a WF $(i,k)\mapsto \vphi (i,k)\geq 0$,
a CF $(i,k)\mapsto\tq_j (i,k)>0$ and a no-ruin parameter value $b>0$. Set:
\beq\label{eqg*}g^*(i,k)=g(i,k)-(1+\rho ),\;\;i,k\in M.\eeq

We make a convention that, at times $n=0,1,\ldots$, a part $C_n$ of the current capital
$Z_n$ is invested in the result of the $(n+1)$st trial. This brings a profit/loss $C_ng(\veps_n,\veps_{n+1})$
at time $n+1$. Thus, we deal with an RF $g(i,k)$ where $i$ is the outcome of the previous trial (which
you know) and $k$ the outcome of the next trial (which is unknown at the time of investment). The
remaining part of the capital, $Z_n-C_n$, is invested in the riskless asset. It means that,
compared with \eqref{eq13B1}, we have
the recursion: for $n\geq 1$,
\beq\label{eqRec221}
Z_n=(1+\rho )(Z_{n-1}-C_{n-1})+C_{n-1}g(\veps_{n-1},\veps_n) =
Z_{n-1}\left[1+\rho+\diy\frac{C_{n-1}}{Z_{n-1}}g^*(\veps_{n-1},\veps_n)\right].\eeq
Next, for given $b>0$ and CFs $(i,k)\mapsto \tq_j(i,k)>0$, we consider strategies/policies $C_j$,
$j\geq 0$, such that
\beq\label{eqa1-3.5}\beal{\rm{(a0)}}\;\; C_j\in\fW_j,\;\; {\rm{(a1)}}\;\; 0\leq C_j<
Z_j,\;{\rm{(a2)}}\;\; 1+\rho +\diy\frac{C_jg^*(\veps_j,\veps_{j+1})}{Z_j}\geq b,\\
\hbox{and }\;{\rm{(a3)}}\quad C_j(\beps_0^j)\sum\limits_{l\in M}\vphi (\veps_j,l)\tq_{j+1}(\veps_j,l)
g^*(\veps_j,l)=0.\ena\eeq

Following the same pattern as before, set:
$$\beal S_n=\sum\limits_{j=1}^n\vphi (\veps_{j-1},\veps_j)\ln{\diy\frac{Z_j}{Z_{j-1}}}
=\sum\limits_{j=1}^n\vphi (\veps_{j-1},\veps_j)\ln\,\left[1+\rho
+\diy\frac{C_{j-1}g^*(\veps_{j-1},\veps_j)}{Z_{j-1}}\right].\ena$$
and
$$\rA_n=\sum\limits_{j=1}^n\alpha_{j-1}(\veps_{j-1})\;\hbox{ where }\alph_{j-1}(\veps_{j-1})
=\bbE\left[\vphi (\veps_{j-1},\veps_j)\ln\,\diy\frac{\tp (\veps_{j-1},\veps_j)}{
\tq_j(\veps_{j-1},\veps_j)}\Big|\fW_{j-1}\right].$$
Again we are interested in maximizing the mean value $\bbE S_n$ in $C_0$, $\ldots$, $C_{n-1}$.
Assume the condition that, $\forall$ \ $i\in M$, $j\geq 0$,
\beq\label{eqq12Rl} {\rm{(q-p)}}\qquad
\sum\limits_{l\in M} \vphi (i,l)\big[(1+\rho)\tq_j(i,l)-\tp (i,l)\big]\leq 0.\eeq
Assumptions \eqref{eqa1-3.5} mimic \eqref{eqa1-3.3} while \eqref{eqq12Rl} mimics
\eqref{eqq11Rl}; terminological parallels are also notable here.
\vskip .5 truecm

{\bf Theorem 2.1.} {\sl Adopt Eqns \eqref{eqg*}--\eqref{eqq12Rl}. Then:

{\rm{(a)}} The sequence $\{S_n-\rA_n,\,n\geq 1\}\;$ is a supermartingale  and hence
$\bbE\,S_n\leq\bbE\rA_n$ \ $\forall$ $n\geq 1$.

{\rm{(b)}} To achieve equality $\bbE\,S_n=\bbE \rA_n$: the sequence
$\{S_n-\rA_n\}$ is a martingale for a sequence of RVs $C_j$
satisfying \eqref{eqa1-3.5} iff the following conditions {\rm{(i)}}--{\rm{(ii)}} are fulfilled.

\quad {\rm{(i)}} There exists a function $D: M\to \bbR$ such that, \ $\forall$ \ $i,k\in M$,
\beq\label{eq8M1n}\beal\quad{\rm{(i1)}}\;\; 0\leq D(i)< 1,\;\; {\rm{(i2)}}\;\;
1+\rho+D(i)g^*(i,k)\geq b,\\
\quad{\rm{(i3)}}\;\;D(i)\sum\limits_{l\in M}
{\diy\frac{\vphi (i,l)\tp (i,l)g^*(i,l)}{1+\rho+D(i)g^*(i,l)}}= 0,\;
\hbox{ i.e., either}\\
\qquad\quad {\rm{(i3A)}}\;\; \sum\limits_{l\in M}{\diy\frac{\vphi (i,l)\tp (i,l)g^*(i,l)}{1+\rho+D(i)g^*(i,l)}}= 0\;
\hbox{ or }\;{\rm{(i3B)}}\;\;D(i)=0,\;\hbox{ and}\\
\quad{\rm{(i4)}}\;\hbox{ the CFs \ $\tq_j$ \ have the form }\;
\tq_j (i,k)=\diy\frac{\tp (i,k)}{1+\rho+D(i)g^*(i,k)},\\
\qquad\quad\hbox{ independently of }\;j\geq 0.\quad\ena\eeq

\quad{\rm{(ii)}} $\forall$ \ $n\geq 1$, \ the policy $C_{n-1}$ produces a proportional investment
portfolio:
after the $(n-1)$st trial the amount $C_{n-1}(\beps_0^{n-1})=D(\veps_{n-1})Z_{n-1}$ goes to
the asset with return $g$ and $\big[1-D(\veps_{n-1})\big]Z_{n-1}$ to the riskless return.
\medskip

{\rm{(c)}} Define the map $i\in M\mapsto D^\rO(i)$ as follows. Given $i$, consider Eqn {\rm{(i3A)}}:
it has at most one solution $D(i)>0$. If {\rm{(i3A)}} has a solution $D(i)>0$ obeying 
conditions {\rm{(i1)}}--{\rm{(i2)}}, set $D^\rO (i)=D(i)$; otherwise
$D^\rO (i)=0$. Then the policy $C^\rO_{n-1}=D^\rO(\veps_{n-1})Z_{n-1}$ yields
the following value $E_n$ for the expectation $\bbE S_n$:
\beq\label{eqOptRaten}\beal E_n=\sum\limits_{j=1}^n\bet_{j-1}\;\hbox{ where }\;
\bet_{j-1}=\bbE\Big\{\vphi (\veps_{j-1},\veps_j)
\ln\,\Big[1+\rho +D(\veps_{j-1})g^*(\veps_{j-1},\veps_j)\Big]\Big\}.\ena\eeq
The maximum of $E_n$ in \eqref{eqOptRaten} over maps $i\in M\mapsto D(i)$
satisfying {\rm{(i1)}}--{\rm{(i3)}} from \eqref{eq8M1n} gives the maximum of $\bbE S_n$ over the
portfolios $\{C_j\}$ satisfying properties {\rm{(a0)}}--{\rm{(a3)}} in Eqn \eqref{eqa1-3.5}.
\medskip

{\rm{(d)}} Suppose that the map $i\in M\mapsto D^\rO(i)$ from assertion {\rm{(c)}} is such that
$D^\rO(i)>0$ (so the alternative {\rm{(i3A)}} holds) $\forall$
$i\in M$. Then the policy $C^\rO_{n-1}=D^\rO(\veps_{n-1})Z_{n-1}$ maximizes
each summand $\bet_{j-1}$ in \eqref {eqOptRate1}, and therefore yields the maximum
of the whole  sum \ $\bbE S_n$, among strategies satisfying properties {\rm{(a0)}}--{\rm{(a2)}}
in \eqref{eqa1-3.5}.}

\vskip .5 truecm

{\bf Proof.} (a) We follow the same pattern as in Theorem 1.1(a). Write:
$$\beal\bbE\Big\{\big(S_n-\rA_n\big)\big|\fW_{n-1}\Big\}=S_{n-1}-\rA_{n-1}\\
\qquad\qquad+\diy\bbE\left\{\left[
\vphi (\veps_{n-1},\veps_n)\ln\,\left(1+\rho+\diy\frac{C_{n-1}g^*(\veps_{n-1},\veps_n)}{Z_{n-1}}
\right)\right]\Big|\fW_{n-1}\right\}-\alph_{n-1}(\veps_{n-1}).\ena$$
Next, represent
\beq\label{eq10*}
\beal\diy\bbE\left\{\left[\vphi (\veps_{n-1}, \veps_n)\ln\,\left(1+\rho+\frac{C_{n-1}g^*(\veps_{n-1}, \veps_n)}{Z_{n-1}}\right)\right]\Big|\fW_{n-1}\right\}-\alph_{n-1} \\
\quad\diy =\;\sum\limits_{l\in M}\vphi_n(\veps_{n-1},l)\;\tp (\veps_{n-1},l)\;\ln\;\left[1+\rho+
\frac{C_{n-1}g^*(\veps_{n-1},l)}{Z_{n-1}}\right]\\
\qquad\qquad\diy - \sum\limits_{l\in M}\vphi (\veps_{n-1},l)\tp (\veps_{n-1},l)\ln\,
\frac{\tp (\veps_{n-1},l)}{\tq_n(\eps_{n-1},l)}\\
\quad\diy = \sum\limits_{l\in M}\vphi(\veps_{n-1},l)\tp (\veps_{n-1},l) \ln\frac{1+\rho
+C_{n-1}g^*(\veps_{n-1},l)/Z_{n-1}}{\tp (\veps_{n-1},l)/\tq_n(\veps_{n-1},l )}\\
\quad\diy = \sum\limits_{l\in M}\vphi (\veps_{n-1},l)\;\tp (\veps_{n-1},l)\;
\ln\;\frac{h_n(\veps_{n-1},l)}{\tp (\veps_{n-1},l)}\\
\quad\diy \leq\sum\limits_{l\in M}\vphi (\veps_{n-1},l)\tp (\veps_{n-1},l)
\left[ \frac{h_n(\veps_{n-1},l)}{\tp (\veps_{n-1},l)} - 1 \right]{\mathbf 1}\big(\tp
(\veps_{n-1},l)>0\big)\\
\quad\diy =\sum\limits_{l\in M}\vphi (\veps_{n-1},l) \Big[h_n(\veps_{n-1},l)
- \tp (\veps_{n-1},l)\Big]\leq 0.\ena\eeq
Here $h_n(\veps_{n-1},k) := \diy\tq_n(\veps_{n-1},k )\left[1 +\rho+\frac{C_{n-1}g^*(\veps_{n-1},k )
}{Z_{n-1}}\right]$, $k\in M$.

The final inequality in \eqref{eq10*} holds since
\beq\label{eq11.1*}\beal
\diy \sum\limits_{l\in M}\vphi (\veps_{n-1},l)h_n(\veps_{n-1},l ) = (1+\rho )\sum\limits_{l\in M}
\vphi (\veps_{n-1},l)\tq_n(\veps_{n-1},l )\\
\qquad+\diy\frac{C_{n-1}}{Z_{n-1}}\sum\limits_{l\in M}\vphi (\veps_{n-1},l)\tq_n(\veps_{n-1},l)
g^*(\veps_{n-1},l)\leq
\sum\limits_{l\in M}\vphi (\veps_{n-1},l)
\tp (\veps_{n-1},l),\ena\eeq
because of properties (a3) in Eqn \eqref{eqa1-3.5} and (q--p) in Eqn \eqref{eqq12Rl}.

(b) As before, the martingale property emerges when we have equality in the inequalities in
\eqref {eq10*}. The analysis of these situations proceeds as in the proof of Theorem 1.1(b).

(c,d). The proof of these assertions does not differ from the proof of their counterparts
in Theorem 1.1. $\Box$

\vskip .5 truecm

{\bf Example 2.1: IID trials.} Consider IID trials in the context of Scheme II.
As in Example 1.1, we now have $\tp (i,k)=\tp(k)$
and work with $g(i,k)=g(k)$, $\vphi (i,k)=\vphi (k)$ and $\tq_j(i,k)=\tq(k)$. Choose
$b\in (0,1+\rho )$. As in \eqref{eqg*}, set: $g^*(k)=g(k)-(1+\rho )$. Eqns \eqref{eqRec221}
and \eqref{eqa1-3.5} read: $\forall$ $n\geq 1$,
\beq\label{eqSust2A}\beal
Z_n=(1+\rho )(Z_{n-1}-C_{n-1})+C_{n-1}g(\veps_n) =
Z_{n-1}\left[1+\rho+\diy\frac{C_{n-1}g^*(\veps_n)}{Z_{n-1}}\right]\;\hbox{ where}\\
{\rm{(a0)}}\;\; C_{n-1}\in\fW_{n-1},\;\; {\rm{(a1)}}\;\; 0\leq C_{n-1}< Z_{n-1},
\;{\rm{(a2)}}\;\; 1+\rho +\diy\frac{C_{n-1}g^*(\veps_n)}{Z_{n-1}}\geq b,\hbox{ and}\\\\
{\rm{(a4)}}\quad C_{n-1}\sum\limits_{l\in M}\vphi (l)\tq (l)g^*(l)=0\;\
\hbox{ (meaning that at least one factor is $0$).}\ena\eeq

We now define:
$$\beal S_n=\sum\limits_{j=1}^n\vphi (\veps_j)\ln\diy\frac{Z_j}{Z_{j-1}}
=\sum\limits_{j=1}^n\vphi (\veps_j)\ln\,\left[1+\rho+\diy\frac{C_{j-1}g^*(\veps_j)}{Z_{j-1}}
\right]\ena$$
and
$$\rA_n=n\alpha\;\hbox{ where }\;\alph =\bbE\left[\vphi (\veps_1)\ln\,\diy\frac{\tp (\veps_1)}{
\tq(\veps_1)}\right].$$
Again we are interested in maximizing the mean value $\bbE S_n$ in $C_0$, $\ldots$, $C_{n-1}$.
Eqns \eqref{eqq12Rl} and \eqref{eq8M1n} take the form, respectively,
\beq\label{eqq12Rn}\beal
{\rm{(q-p)}}\qquad
\sum\limits_{l\in M} \vphi (l)\big[(1+\rho)\tq (l)-\tp (l)\big]\leq 0,\;\hbox{ and}\\
{\rm{(i1)}}\;\;\exists\;\;\hbox{ a constant }\;D\in [0,1)\;\hbox{such  that}\;{\rm{(i2)}}\;\;
1+\rho+Dg^*(k)\geq b,\\
{\rm{(i3)}}\;\;D\sum\limits_{l\in M}{\diy\frac{\vphi (l)\tp (l)g^*(l)}{1+\rho+Dg^*(l)}}= 0,\;
\;{\rm{(i4)}}\;\hbox{the CF is}\;\tq (k)=\diy\frac{\tp (k)}{1+\rho+Dg^*(k)},\;k\in M.\ena\eeq
The expression for $\bbE S_n$ in Eqn \eqref{eqOptRaten} is
$$\beal E_n=n\bet\;\;\hbox{ where }\;\;\bet =\bbE\Big\{\vphi (\veps_1)
\ln\,\Big[1+\rho +Dg^*(\veps_1)\Big]\Big\},\ena$$
which should be maximized over all choices of $D$ satisfying (i1)--(i3).

As before, to make the conditions more explicit, we take $M=\{0,1\}$. To simplify,
set $g(1)=-g(0)=\gam >0$ and $\vphi (1)=\vphi (0)=1$, with $g^*(1)=\gam -(1+\rho )$,
$g^*(0)=-\gam -(1+\rho )$. Take $0<b\leq 1+\rho$. Re-write Eqn \eqref{eqq12Rn}:
\beq\label{eqspec4}\beal{\rm{(q-p)}}\;\;(1+\rho )\big[\tq (1)+\tq (0)\big]\leq 1;\\
{\rm{(i1)}}\;\;\exists\;D\in [0,1)\;\hbox{ such that }\;\;{\rm{(i2)}}\;\;
D\leq\diy\frac{1+\rho-b}{\gam +1+\rho},\\
{\rm{(i3)}}\;D\left[
{\diy\frac{\tp (1)(\gam -1-\rho )}{1+\rho+D(\gam -1-\rho )}}
-{\diy\frac{\tp (0)(\gam +1+\rho )}{1+\rho-D(\gam +1+\rho )}}\right]=0,\;\hbox{ and}\\
{\rm{(i4)}}\;\;\hbox{the CF values are }\;
\tq (1)=\diy\frac{\tp (1)}{1+\rho+D(\gam -1-\rho )},\;
\tq (0)=\diy\frac{\tp (0)}{1+\rho-D(\gam +1+\rho )}.\quad\ena\eeq

Now, Eqn (i3) in \eqref{eqspec4} can be solved explicitly:
% $$\tp_1(\gam -1-\rho )[1+\rho+D(-\gam -1-\rho )]+\tp_2(-\gam -1-\rho) [1+\rho+D(\gam -1-\rho )]=0$$
% whence
\beq\label{eqDwrho}D=0\;\hbox{ or }\;D=D_0\;\hbox{ where }\;D_0=(1+\rho )\frac{\gam\big[\tp (1)-\tp (0)\big]-(1+\rho)}{\gam^2-(1+\rho)^2}.\eeq
We see that the non-trivial solution $D>0$ emerges iff either (A) $\gam\big[\tp (1)-\tp (0)\big]>1+\rho$ \ or
(B) $\gam <1+\rho$.  Next, taking $\tq$ as in (i4) yields equality in (q--p) which is equivalent to (i3).

Thus, it remains to check the relations (i1)--(i2) for $D$ from \eqref{eqDwrho}:
\beq\label{eqNonTD1}0\leq (1+\rho )\frac{\gam\big[\tp (1)-\tp (0)\big]-(1+\rho)}{\gam^2-(1+\rho)^2}\leq 1\wedge
\frac{1+\rho-b}{\gam +1+\rho}.\eeq
If, for given $\tp (1), \tp (0)$, $\gam$ and $\rho$, the inequalities in \eqref{eqNonTD1} are satisfied
then we have two solutions for $D$ specified in \eqref{eqDwrho}. Next, we have to compare
$\bet =\bbE\vphi (\veps_1)\ln\,\big[1+\rho +D_0g^*(\veps_1)\big]$ for $D=D_0$ and
$\bet =\ln\,(1+\rho )\bbE\vphi (\veps_1)$ for $D=0$: the larger solution
identifies the optimizer $D^\rO$ giving a maximal value for $\bbE S_n$.  Viz., such a
comparison shows that $D^\rO=0$ in the situation (B) \ $\gamma<1+\rho$, in which case the
investor allocates all capital to the riskless asset.

\section{Markov trials with several risky assets}

Now assume that, within the remits of the investment Scheme I, we have to make a choice between $K$
risky assets, with a vector of individual RFs $\ug =\left(g^{(1)},\ldots ,g^{(K)}\right)$. As before, let the trial
results $\veps_n$ be generated by an MC with states $i,k$ from space $M$, finite or countable, and with
transition probabilities
$\tp (i,k)$. As before, the RF value $g^{(s)}(i,k)$ gives the return from asset $s$ at the time
(or immediately after) the $n$th trial when the outcome is $k$ preceded by outcome $i$ at the $(n-1)$st trial.
Here we consider a sequence  $\{\uC_n,\,n\geq 0\}$ where $\uC_n=\left(C^{(1)}_n,\ldots ,C^{(K)}_n\right)$ is an
$\fW_n$-measurable random $K$-dimensional vector (a predictable policy/strategy portfolio). The recursion
for $Z_n$ is similar to \eqref{eq13B1}, with replacing scalar random variables by random vectors (RVs):
\beq\label{eq13B2}
Z_n=Z_{n-1}+\uC_{n-1}\ug(\veps_{n-1},\veps_n)=Z_{n-1}\left[1
+\frac{\uC_{n-1}\ug(\veps_{n-1},\veps_n)}{Z_{n-1}}\right],\;n\geq 1.
\eeq
Here and below,
$$\uC_{n-1}\ug (\veps_{n-1},\veps_n)=\sum\limits_{s=1}^KC^{(s)}_{n-1}g^{(s)}(\veps_{n-1},\veps_n).$$
Also, $|\uC_j|=\sum\limits_{s=1}^KC^{(s)}_j=\uC_{n-1}\u1$ where $\u1=(1,\ldots ,1)$, and we write
$\uC_j\geq\U0$ if $C^{(s)}_j\geq 0$ \ $\forall$ \ $s$. A similar convention is in place for other expressions
of an analogous structure.

% Next, introduce a weight function, $(i,k)\mapsto \vphi (i,k)\in\bbR$,
% representing a utility value assigned to outcome $k$ when it succeeds a previous outcome $i$.
% \newpage

As above, we set $S_n:=\sum\limits_{j=1}^n\vphi (\veps_{j-1},\veps_j)\ln\,\diy\frac{Z_j}{Z_{j-1}}$
and aim at maximizing the mean value $\bbE S_n$ in $\uC_0,\ldots ,\uC_{n-1}$, under certain
restrictions. Fix $b>0$ and re-write the conditions outlining the portfolio classes under consideration:
\ $\forall$ \ $j\geq 0$,
\beq\label{eqa1-3.6}\beal{\rm{(a0)}}\;\uC_j\in\fW_j\;\hbox{(predictability)},\; {\rm{(a1)}}\; \uC_j\geq\U0,\;|\uC_j|<
Z_j\;\hbox{(sustainability)},\\
{\rm{(a2)}}\;1+\diy\frac{\uC_j\ug(\veps_j,\veps_{j+1})}{Z_j}\geq b\;\hbox{ (no ruin), \ and}\\
{\rm{(a3)}}\quad \sum\limits_{l\in M}\vphi (\veps_j,l)\tq_{j+1}(\veps_j,l)\left[\uC_j(\beps_0^j)
\ug(\veps_j,l)\right]=0\;\hbox{ (weighted (q,g)-balance).}\ena\eeq
% \;\hbox{ which can be related to}\\
% {\rm{(a4)}}\;\sum\limits_{l\in M}\vphi (i,l)\tq_n(i,l)g(i,l)=0,\;\hbox{ and }\;{\rm{(a5)}}\;C_{n-1}
% (\bi_0^{n-1})=0,\;\hbox{ for some or all }\\
% \qquad i\in M,\;\bi_0^{n-1}=(i_0,\ldots ,i_{n-1})\in M^n.\ena\eeq
We also assume, \ $\forall$ \ $i\in M$, the (q--p) bound \eqref{eqq11Rl}.
%\beq\label{eqq11Rl} {\rm{(q-p)}}\qquad
%\sum\limits_{l\in M} \vphi (i,l)\big[\tq_j (i,l)-\tp (i,l)\big]\leq 0.\eeq

We again want to determine, when possible, a sequence of optimal strategies. For
$\alph_j(i)$ and $\rA_n$ we
follow the definitions from Eqns \eqref{eq17B1} and \eqref{eq101M}.

It is instructive to observe that the vector case can be treated as scalar after we fix
the fractions $h_j^{(s)}:=C_j^{(s)}\big/|\uC_j|$, $1\leq s\leq K$, and introduce the
`weighted' RFs: ${\ov g}(i,k)=\sum\limits_{s=1}^K
h^{(s)}_jg^{(s)}(i,k)$. Such a view will be useful when one considers examples; see below.

Consider the following conditions (i1)--(i4) which are vector counterparts of their
scalar predecessors from \eqref{eq8M1}. For convenience, we use the same labelling system
as above.

\quad {\rm{(i)}}\; There exists a map $i\in M\mapsto \uD (i)$ where vector
$\uD(i)=\left(D^{(1)}(i),\ldots ,D^{(K)}(i)\right)$ is such that \ $\forall$ \ $i,k\in M$,
\beq\label{eq8M133}\beal\quad{\rm{(i1)}}\quad \uD (i)\geq\U0,\;\hbox{ and }\;\;|\uD(i)|<1\;\;
\hbox{ (D-sustainability)},\\
\quad{\rm{(i2)}}\quad  1+\uD (i)\cdot\ug(i,k)\geq b\;\hbox{ (D-no-ruin),}\\
\quad{\rm{(i3)}}\quad\sum\limits_{l\in M}\vphi (i,l)\tp (i,l)\,\diy\frac{\uD (i)\cdot\ug (i,l)}{1
+\uD (i)\cdot\ug (i,l)}
= 0\;\hbox{(WE (D,g)-balance), and}\\
\quad{\rm{(i4)}}\quad\;\hbox{the CFs \ $\tq_j$ \ are }\;\; \tq_j(i,k)=\diy\frac{\tp (i,k)}{1+\uD (i)\cdot\ug(i,k)},\\
\qquad\qquad\hbox{independently of }\;j\geq 0\;\hbox{ (q-representation).}\ena\eeq
Here the analog of the first alternative in (i3) is that $\uD(i)$ satisfies a system of equations: \
\beq\label{eqM8133}\beac
\quad\;\;{\rm{(i3A)}}\quad\sum\limits_{l\in M}\tp (i,l)\,\diy\frac{\vphi (i,l)g^{(s)}(i,l)}{1+\uD (i)\cdot\ug (i,l)}=0\qquad\qquad\qquad\qquad\qquad\qquad\qquad\qquad{}\\
\qquad\qquad\qquad\qquad\hbox{$\forall$ \
$i\in M$ and $1\leq s\leq K$ (strong WE (D,g)-balance).\qquad}\ena\eeq
Cf. (i3A) in Eqn \eqref{eq8M1}.

\vskip .5 truecm

{\bf Theorem 3.1.} {\sl Assume the above setting \eqref{eq13B2}--\eqref{eqM8133}. The following
assertions hold true.

{\rm{(a)}} Take any sequence $\{\uC_j,\,j\geq 0\}$ obeying {\rm{(a0)}}--{\rm{(a3)}} in \eqref{eqa1-3.6}.
Then the sequence $\{S_n-\rA_n,\,n\geq 1\}\;$ is a supermartingale; hence $\bbE\,S_n\leq \,
\bbE \rA_n$ \ $\forall$ \ $n\geq 1$.

{\rm{(b)}} To reach equality $\bbE\,S_n=\sum\limits_{j=1}^n\bbE \alph (\veps_{j-1})$: the sequence
$\{S_n-\rA_n\}$ is a martingale for a sequence of RVs \ $\{\uC_j\}$, satisfying {\rm{(a0)}}--{\rm{(a3)}} \ iff \
the additional properties {\rm{(i)}}, {\rm{(ii)}} below are fulfilled. \vskip .3 truecm

\quad {\rm{(i)}}  There exists a map $i\in M\mapsto \uD (i)$ where vector
$\uD(i)=\left(D^{(1)}(i),\ldots ,D^{(K)}(i)\right)$ is such that \ $\forall$ \ $i,k\in M$, properties
{\rm{(i1)}}--{\rm{(i3)}} in \eqref{eq8M133} are fulfilled, and the CFs $\tq_j$ are as in {\rm{(i4)}}.

\quad{\rm{(ii)}} The portfolio vectors $\uC_j$ have components $C^{(s)}_j(\beps_0^j)=D^{(s)} (\veps_j)Z_j$,
$1\leq s\leq K$, $j\geq 0$. That is, the prescribed fractions of the capital value $Z_j$
are invested in the available returns.

{\rm{(c)}} Suppose there exists a map $i\in M\mapsto\uD(i)=\left(D^{(1)}(i),\ldots ,D^{(K)}(i)\right)$
fulfilling conditions {\rm{(i1)}}--{\rm{(i3)}} in Eqn \eqref{eq8M133}. Let \ $\bD$ \ stand for the array
of values $D^{(s)}(i)$, $i\in M$, $1\leq s\leq K$, and define the quantity $E_n=E_n(\bD)$ by
\beq\label{eqOptRateM}\beal E_n=\sum\limits_{j=1}^n\bet_{j-1}\;\hbox{ where }\;
\bet_{j-1}=\bbE\Big\{\vphi (\veps_{j-1},\veps_j)
\ln\,\Big[1+\uD(\veps_{j-1})\cdot\ug(\veps_{j-1},\veps_j)\Big]\Big\}.\ena\eeq
Consider the optimization problem
\beq\label{eqmAxiM1}{\rm{maximise}}\quad E_n(\bD )\quad\hbox{subject to } {\rm{(i1)-(i3)}.}\eeq
Let \ $\bD^*={\rm{arg\;max}}\;E_n$ be a (possibly, non-unique) optimizer,  and $E^*_n=E_n(\bD^*)$
denote the optimal value for \eqref{eqmAxiM1}. Then $E^*_n$ defines the maximum of the
expectation \ $\bbE S_n$ \ among all portfolios $\{\uC_j\}$ satisfying the
properties {\rm{(a0)}}--{\rm{(a3)}} in Eqn \eqref{eqa1-3.6}.  The optimizer \ $\bD^*$ \ written as
a collection of vectors $\uD^*(i)$, $i\in M$, yields a proportional investment portfolio where $\uC_j(\beps_0^j)=
\uD^*(\veps_j)Z_j(\beps_j^j)$.

{\rm{(d)}} Suppose there exists a map $i\in M\mapsto D(i)$ fulfilling conditions {\rm{(i1)}}--{\rm{(i2)}}
and {\rm{(i3A)}} in Eqns \eqref{eq8M133} and \eqref{eqM8133}, respectively. Then such a map is unique,
and the proportional investment portfolio $\uC^\rO_{j-1}=\uD^{*}(\veps_{j-1})Z_{j-1}$ maximises
each summand $\bet_{j-1}$ in \eqref {eqOptRateM}. Therefore, it yields the maximum
of the whole  sum \ $\bbE S_n$, among strategies satisfying properties {\rm{(a0)}}--{\rm{(a2)}}
in Eqn \eqref{eqa1-3.6}.}

\medskip

{\bf Proof.} (a) We just repeat the argument from the proof of Theorem 1.1(a) for the vector case. Here
$$\beal\bbE\Big\{\big(S_n-\rA_n\big)\big|\fW_{n-1}\Big\}=S_{n-1}-\rA_{n-1}\\
\qquad +\diy\bbE\left\{\left[
\vphi (\veps_{n-1},\veps_n)\ln\,\left(1+\diy\frac{\uC_{n-1}\ug(\veps_{n-1},\veps_n)}{Z_{n-1}}
\right)\right]\Big|\fW_{n-1}\right\}-\alph_{n-1}(\veps_{n-1}),\ena$$
and we write
\beq\label{eq11*}
\beal\diy\bbE\left\{\left[\vphi (\veps_{n-1}, \veps_n)\ln\,\left(1+\frac{\uC_{n-1}\ug
(\veps_{n-1}, \veps_n)}{Z_{n-1}}\right)\right]\Big|\fW_{n-1}\right\}-\alph_{n-1} \\
\quad\diy =\;\sum\limits_{l\in M}\vphi_n(\veps_{n-1}, l)\;\tp (\veps_{n-1},l)\;\ln\;\left[1+
\frac{\uC_{n-1}\ug (\veps_{n-1},l)}{Z_{n-1}}\right]\\
\qquad\qquad\diy - \sum\limits_{l\in M}\vphi (\veps_{n-1},l)\tp (\veps_{n-1},l)\ln\,
\frac{\tp (\veps_{n-1},l)}{\tq_n (\eps_{n-1},l)}\\
\quad\diy = \sum\limits_{l\in M}\vphi(\veps_{n-1},l)\tp (\veps_{n-1},l) \ln\frac{1
+\uC_{n-1}\ug (\veps_{n-1},l)/Z_{n-1}}{\tp (\veps_{n-1},l)/\tq_n (\veps_{n-1},l )}\\
\quad\diy = \sum\limits_{l\in M}\vphi (\veps_{n-1},l)\;\tp (\veps_{n-1},l)\;
\ln\;\frac{h (\veps_{n-1},l)}{\tp (\veps_{n-1},l)}\\
\qquad\diy \leq\sum\limits_{l\in M}\vphi (\veps_{n-1},l)\tp (\veps_{n-1},l)
\left[ \frac{h_n(\veps_{n-1},l)}{\tp (\veps_{n-1},l)} - 1 \right]{\mathbf 1}\big(\tp
(\veps_{n-1},l)>0\big)\\
\quad\diy =\sum\limits_{l\in M}\vphi (\veps_{n-1},l) \Big[h_n(\veps_{n-1},l)
- \tp (\veps_{n-1},l)\Big]\leq 0.\ena\eeq
Here $h_n(\veps_{n-1},k) := \diy\tq_n(\veps_{n-1},k )\left[1
+\frac{\uC_{n-1}\ug (\veps_{n-1},k )
}{Z_{n-1}}\right]$, $k\in M$. The final inequality in \eqref{eq11*} holds since
\beq\label{eq12.1*}\beal
\diy \sum\limits_{l\in M}\vphi (\veps_{n-1},l)h_n(\veps_{n-1},l ) = \sum\limits_{l\in M}
\vphi (\veps_{n-1},l)\tq_n (\veps_{n-1},l )\\
\qquad+\diy\frac{1}{Z_{n-1}}\sum\limits_{l\in M}\vphi (\veps_{n-1},l)\tq_n(\veps_{n-1}, l)
\uC_{n-1}\cdot\ug(\veps_{n-1},l)\leq
\sum\limits_{l\in M}\vphi (\veps_{n-1},l)\tp (\veps_{n-1},l),\ena\eeq
because of property (q--p) in Eqn \eqref{eqq11Rl}.

(b) The proof of this assertion is reduced to the analysis of the equality cases in
\eqref{eq11*}. It  does not differ from assertion (b) in Theorem 1.1.

The proof of assertion (c) is again a straightforward inspection.

For (d), we can repeat the argument used in the proof of Theorem 1.1(d).
It yields only some notational complications without making the situation
different in principle.  $\Box$

\vskip .5 truecm

Theorems 2.1 and 3.1
can be combined into a statement about Scheme II, with one riskless and several risky assets.
The assumption is that the trials results $\veps_j$, $j\geq 0$, are still generated by an MC
on $M$ with transition probabilities $\tp (i,k)$, and we have $K+1$ assets with RFs
$g^{(0)}$ and $\ug=\left(g^{(1)},\ldots ,g^{(K)}\right)$. The asset with RF $g^{(0)}$ is
riskless: $g^{(0)}(i,k)=1+\rho$ with $\rho >0$. The assets with RFs $g^{(1)},\ldots ,g^{(K)}$
are risky: they can bring profit or loss $g^{(s)}(i,k)$, $i,k\in M$, $1\leq s\leq K$. We again
adopt the Scheme II: the capital not invested in the risky assets is put in the riskless return.

Observe that the model under consideration emerges as a special case of the
previous model if we include $g^{(0)}$ into vector $\ug$, increasing its dimension from $K$
to $K+1$. However, the explicit use of the form $\rho^{(0)}(k)=1+\rho$ of the
riskless asset makes the presentation less abstract.

As before, let $\uC_{n-1}$ denote the random vector $\left(C^{(1)}_{n-1},\ldots ,C^{(K)}_{n-1}\right)$
where $C^{(s)}_{n-1}$ stands for the amount of investment in the $s$th risky asset before
the $n$th trial, $n\geq 1$.
% As usually, we will assume predictability: $\uC_{n-1}\in\fW_{n-1}$. Set:
% $|\uC_{n-1}|=\sum\limits_{p=1}^sC^{(p)}_{n-1}$.
The recursion for the capital value becomes
\beq\label{eqRL&sR}Z_n=(1+\rho )(Z_{n-1}-|C_{n-1}|)+\uC_{n-1}\ug(\eps_{n-1},\eps_n)=Z_{n-1}\left[
1+\rho +\diy\frac{\uC_{n-1}\ug^*(\veps_{n-1},\veps_n)}{Z_{n-1}}\right].\eeq
Here and below in this section:
\beq\label{eqg*gen}\ug^*(i,k)=\left(g^{(1)}(i,k)-(1+\rho ),\ldots ,g^{(K)}(i,k)-(1+\rho )\right)
=\ug (i,k)-(1+\rho )\u1.\eeq
Cf. Eqns \eqref{eqg*}, \eqref{eqRec221} and \eqref{eq13B2}.

The assumption on the sequence of portfolios $\{\uC_j\}$ and functions $\ug$, $\vphi$ and
$\tq_j$ are: \ $\forall$ \ $j\geq 0$,
\beq\label{eqa1-3.7}\beal{\rm{(a0)}}\;\uC_j\in\fW_j,\;\; {\rm{(a1)}}\; \uC_j\geq\U0,\;|\uC_j|<
Z_j,\;{\rm{(a2)}}\;1+\rho+\diy\frac{\uC_j\ug^*(\veps_j,\veps_{j+1})}{Z_j}\geq b,\\
\hbox{and }\;\;
{\rm{(a3)}}\quad \sum\limits_{l\in M}\vphi (\veps_j,l)\tq_{j+1}(\veps_j,l)\left[\uC_j(\beps_0^j)
\ug^*(\veps_j,l)\right]=0,\ena\eeq
again citing sustainability, no-ruin and balance conditions. The dominance condition \eqref{eqq12Rl}
will also be used.

The quantity of interest is still the mean value $\bbE S_n$. The goal is to maximise
$\bbE S_n$ in $\uC_1,\ldots ,\uC_{n-1}$ subject to restrictions \eqref{eqa1-3.7}, with or
without condition (a3).

To this end, we list the adapted conditions (i1)--(i4):

\quad {\rm{(i)}}\; There exists a map $i\in M\mapsto \uD (i)$ where vector
$\uD(i)=\left(D^{(1)}(i),\ldots ,D^{(K)}(i)\right)$ is such that \ $\forall$ \ $i,k\in M$,
\beq\label{eq8M135}\beal\quad{\rm{(i1)}}\quad \uD (i)\geq\U0,\;\hbox{ and }\;\;|\uD(i)|<1\;\;
\hbox{ (D-sustainability)},\\
\quad{\rm{(i2)}}\quad  1+\rho +\uD (i)\cdot\ug^*(i,k)\geq b\;\hbox{ (D-no-ruin),}\\
\quad{\rm{(i3)}}\quad\sum\limits_{l\in M}\vphi (i,l)\tp (i,l)\,\diy\frac{\uD (i)\cdot\ug^*(i,l)}{1+\rho
+\uD (i)\cdot\ug^*(i,l)}
= 0\;\hbox{(WE (D,g)-balance), and}\\
\quad{\rm{(i4)}}\quad\;\hbox{the CFs \ $\tq_j$ \ are }\;\; \tq_j(i,k)=\diy\frac{\tp (i,k)}{1+\rho+\uD (i)\cdot\ug^*(i,k)},\\
\qquad\qquad\hbox{independently of }\;j\geq 0\;\hbox{ ($\tq$-representation),}\ena\eeq
with the alternative in (i3A) as follows:
\beq\label{eqM8136}\beac
\quad\;\;{\rm{(i3A)}}\quad\sum\limits_{l\in M}\tp (i,l)\,\diy\frac{\vphi (i,l)g^{(s)}(i,l)}{1+\rho +\uD (i)\cdot\ug^* (i,l)}=0\qquad\qquad\qquad\qquad\qquad\qquad\qquad\qquad{}\\
\qquad\qquad\qquad\qquad\hbox{$\forall$ \
$i\in M$ and $1\leq s\leq K$ (strong WE (D,g)-balance).\qquad}\ena\eeq
Cf. (i3A) in Eqn \eqref{eq8M133}.

The cumulative weighted KL entropy process $\{\rA_n\}$ is defined as in Eqns \eqref{eq17B1}
and \eqref{eq101M}.

The combined statement is Theorem 3.2 below. We omit its proof as it repeats
that of Theorems 2.1 and 3.1.

\vskip .5 truecm

{\bf Theorem 3.2.} {\sl Assuming the above setting, we obtain the following assertions.

{\rm{(a)}} Take any sequence $\{\uC_j,\,j\geq 0\}$ obeying {\rm{(a0)}}--{\rm{(a3)}} in \eqref{eqa1-3.7}.
Then the sequence $\{S_n-\rA_n,\,n\geq 1\}\;$ is a supermartingale; hence $\bbE\,S_n\leq \,
\bbE \rA_n$ \ $\forall$ \ $n\geq 1$.

{\rm{(b)}} To reach equality $\bbE\,S_n=\sum\limits_{j=1}^n\bbE \alph (\veps_{j-1})$: the sequence
$\{S_n-\rA_n\}$ is a martingale for a sequence of RVs \ $\{\uC_j\}$, satisfying {\rm{(a0)}}--{\rm{(a3)}} \ iff \
the additional properties {\rm{(i)}}, {\rm{(ii)}} below are fulfilled. \vskip .3 truecm

\quad {\rm{(i)}}  There exists a map $i\in M\mapsto \uD (i)$ where vector
$\uD(i)=\left(D^{(1)}(i),\ldots ,D^{(K)}(i)\right)$ is such that \ $\forall$ \ $i,k\in M$, properties
{\rm{(i1)}}--{\rm{(i3)}} in \eqref{eq8M135} are fulfilled, and the CFs $\tq_j$ are as in {\rm{(i4)}}.

\quad{\rm{(ii)}} The portfolio vectors $\uC_j$ have components $C^{(s)}_j(\beps_0^j)=D^{(s)} (\veps_j)Z_j$,
$1\leq s\leq K$, $j\geq 0$. That is, the prescribed fractions of the capital value $Z_j$
are invested in the available returns while the rest is put in the riskless asset.

{\rm{(c)}} Suppose there exists a map $i\in M\mapsto\uD(i)=\left(D^{(1)}(i),\ldots ,D^{(K)}(i)\right)$
fulfilling conditions {\rm{(i1)}}--{\rm{(i3)}} in Eqn \eqref{eq8M135}. As before, let \ $\bD$ \ stand for the array
of values $D^{(s)}(i)$, $i\in M$, $1\leq s\leq K$, and define the quantity $E_n=E_n(\bD )$ by
\beq\label{eqOptRatem}\beal E_n=\sum\limits_{j=1}^n\bet_{j-1}\;\hbox{ where }\;
\bet_{j-1}=\bbE\Big\{\vphi (\veps_{j-1},\veps_j)
\ln\,\Big[1+\rho +\uD(\veps_{j-1})\cdot\ug^*(\veps_{j-1},\veps_j)\Big]\Big\}.\ena\eeq
Consider the optimization problem
\beq\label{eqmAxiM}{\rm{maximise}}\quad E_n(\bD )\quad\hbox{subject to } {\rm{(i1)-(i3)}.}\eeq
Take an optimizer \ $\bD^*={\rm{arg\;max}}\;E_n$ (possibly, non-unique) and let $E^*_n=E_n(\bD^*)$
denote the optimal value for \eqref{eqmAxiM}.
Then $E^*_n$ defines the maximum of the expectation \ $\bbE S_n$ \ among all portfolios $\{\uC_j\}$ satisfying the
properties {\rm{(a0)}}--{\rm{(a3)}} in Eqn \eqref{eqa1-3.7}.  The array \ $\bD^*$ \ written as
a collection of vectors $\uD^*(i)$, $i\in M$, yields a proportional investment portfolio with
$\uC_j(\beps_0^j)= \uD^*(\veps_j)Z_j(\beps_0^j)$.

{\rm{(d)}} Suppose there exists a map $i\in M\mapsto D(i)>0$ fulfilling conditions {\rm{(i1)}}--{\rm{(i2)}}
and {\rm{(i3A)}} in Eqns \eqref{eq8M135} and \eqref{eqM8136}, respectively. Then such a map is unique,
and the proportional investment portfolio $\uC^\rO_{j-1}=\uD^*(\veps_{j-1})Z_{j-1}$ maximizes
each summand $\bet_{j-1}$ in \eqref {eqOptRatem}. Therefore, it yields the maximum
of the whole  sum \ $\bbE S_n$, among strategies satisfying properties {\rm{(a0)}}--{\rm{(a2)}}
in Eqn \eqref{eqa1-3.7}.}

\vskip .5 truecm

{\bf Example 3.1: IID trials with two risky assets.}  In this example we adopt Scheme I. As was noted,
for IID trials, $\tp (i,k)=\tp (k)$. We again take
$\vphi (i,k)=\vphi (k)$, $\tq_j(i,k)=\tq (k)$ and set \ $\ug(i,k)=\ug (k)$ where $\ug (k)=\left(g^{(1)}(k),g^{(2)}(k)
\right)$. Recall, the portfolio
has the form $\uC_n=(C^{(1)}_n,C^{(2)}_n)$. Conditions  \eqref{eqq11Rl}, \eqref{eqa1-3.6} and
\eqref{eq8M133} are summarized as
\beq\label{eqaVsyo}\beal{\rm{(q-p)}}\quad
\sum\limits_{l\in M} \vphi (l)\big[\tq (l)-\tp (l)\big]\leq 0;\\
{\rm{(a0)}}\;\uC_j\in\fW_j,\;{\rm{(a1)}}\; \uC_j\geq\U0,\;|\uC_j|<
Z_j,\;{\rm{(a2)}}\;1+\diy\frac{\uC_j\cdot\ug(k)}{Z_j}\geq b,\\
\hbox{and }\;{\rm{(a3)}}\; \sum\limits_{l\in M}\vphi (l)\tq (l)\left[\uC_j\cdot
\ug(l)\right]=0,\;\;\forall\;k\in M,\;j\geq 0;\\
{\rm{(i)}}\;\hbox{$\exists$ a vector $\uD =\left(D^{(1)},D^{(2)}\right)$ such that,
\ $\forall$ $k\in M$,}\\
\quad{\rm{(i1)}}\quad \uD\geq\U0\;\hbox{ and }\;\;|\uD|<1,\;\;{\rm{(i2)}}\;\;1+\uD\cdot\ug(k)\geq b,\\
\quad{\rm{(i3)}}\;\sum\limits_{l\in M}\tp (l)\,\diy\frac{\vphi (l)\uD\cdot\ug (l)}{1
+\uD\cdot\ug (l)}
= 0\;\hbox{ and }\;{\rm{(i4)}}\;\hbox{the CF is }\; \tq (k)=\diy\frac{\tp (k)}{1+\uD\cdot\ug(k)},\\
\quad\hbox{with }\;{\rm{(i3A)}}\;\;\sum\limits_{l\in M}\tp (l)\,\diy\frac{\vphi (l)g^{(s)}(l)}{1+\uD\cdot\ug (l)}=0,
\;\;s=1,2.\ena\eeq

In accordance with \eqref{eqmAxiM1}, the optimal portfolio is connected with the maximization problem:
\beq\label{eq8B2}\beal\max\quad\bet (\uD):=\sum\limits_{l\in M}\vphi (l)
\tp (l)\ln\,[1+\uD\cdot\ug(l)]\\
\hbox{in \ $\uD=\left(D^{(1)},D^{(2)}\right)$, subject to (i1)--(i3)}.\ena\eeq
Let $\uD^\rO$ stand for the optimizer: $\uD^\rO={\rm{argmax}}\;\;\bet (\uD)$ in \eqref{eq8B2}.
Note that the Hessian $2\times 2$-matrix $H(\bet )$ is non-positive definite
as its determinant is $0$ and the trace is negative:
$$\beacl H(\bet )&=\Big(\diy\frac{\partial^2\bet (\uD )}{
\partial D^{(s)}\partial D^{(s')}},\;s,s'=1,2\Big)\\
\;&=-\sum\limits_{l\in M}
{\diy\frac{\vphi (l)\tp (l)}{[1+\uD\cdot\ug(l)]^2}}\;\begin{pmatrix}\left[g^{(1)}(l)\right]^2&g^{(1)}(l)g^{(2)}(l)\\
g^{(1)}(l)g^{(2)}(l)&\left[g^{(2)}(l)\right]^2\end{pmatrix}.\ena$$

It shows that $\uD\mapsto \bet (\uD)$ is
a concave function over the polygon $\bbD$ extracted by conditions (i1) and (i2):
$$\beal\bbD:=\Big\{\uD=(D^{(1)},D^{(2)}):\;D^{(1)},D^{(2)}\geq 0,\;D^{(1)}+D^{(2)}<1,\;
1+\uD\cdot\ug(k)\geq b\;\forall\;k\in M\Big\}.\ena$$
We see that if $H(\bet )$ is strictly negative definite then $\uD\in\bbD\mapsto\alph (\uD)$ attains
a maximum at a single point, either in the interior of \ $\bbD$ \
or on the boundary $\partial\bbD$. The stationary point $\uD_{\,0}$, with
$$\nabla\bet (\uD_0)=
\sum\limits_{l\in M}{\diy\frac{\vphi (l)\tp (l)\ug(l)}{1+\uD_0\cdot\ug(l)}}=\U0,$$
is the first candidate
for $\uD^\rO$. Namely, if $\uD_{\,0}\in{\rm{Int}}\;\bbD$ then $\uD^\rO=\uD_0$ as $\uD_0$
solves the equation (i3). In this case, $n\bet (\uD_0)$ yields the maximal value for $\bbE S_n$ among
policies $\{C_j\}$ satisfying (a0)--(a2) in \eqref{eqaVsyo}.

In the case where the stationary point $\uD_0\in\partial\bbD$, the situation is more involved, and we
will not go here into a further detail for a general state space $M$.

However, the corner points $\uD=(0,0)=\U0$, $\uD=(1,0)$ and $\uD=(0,1)$ can be (relatively) easily
checked for optimality (provided that they are feasible for problem \eqref{eq8B2}). Viz., suppose that
$\nabla\bet (\uD)\Big|_{\uD =\U0}=
\sum\limits_{l\in M}\vphi (l)\tp (l)\ug(l)=\U0$.
Then the optimal proportion vector $\uD^\rO$ is $(0,0)$ (no investment by the risk-averse trader). On
the other hand, if $\sum\limits_{l\in M}\vphi (l)g^{(1)}(l)\tp_l>0$ or $\sum\limits_{l\in M}\vphi (l)g^{(2)}(l)\tp_l>0$
then we can look at the gradient values $\nabla\bet (\uD)$ at $\uD=(1,0)$ and $\uD=(1,0)$
(investments into a single asset).

\vskip .5 truecm

To illustrate further, consider again the case $M=\{0,1\}$, with outcomes $k=0,1$ and two probabilities
$\tp (0)$, $\tp (1)$. Without loss of generality, assume $\tp (1)\geq\tp (0)$. Next, for sake of simplicity,
let us again take $\vphi (0)=\vphi (1)=1$. Further, set
$g^{(1)}(1)=-g^{(1)}(0)=\gam_1>0$ and $g^{(2)}(0)=-g^{(2)}(1)=\gam_2>0$, and fix $b\in (0,1)$.

In line with Theorem 3.1(c,d), we seek to solve the optimization problem
\beq\label{eqOPT1}\beal\max\quad\bet (\uD ):=\tp (1)\ln\,\big[1+D^{(1)}\gam_1 -D^{(2)}\gam_2\big]
+\tp (0)\ln\,\big[1-D^{(1)}\gam_1 +D^{(2)}\gam_2\big]\\
\hbox{subject to }\;\uD=(D^{(1)},D^{(2)})\;\hbox{ satisfying }\;D^{(1)},D^{(2)}\geq 0,\;D^{(1)}+D^{(2)}< 1,\\
\qquad\qquad\quad 1+D^{(1)}\gam_1-D^{(2)}\gam_2\geq b,\;1-D^{(1)}\gam_1+D^{(2)}\gam_2\geq b,\;
\hbox{ and}\\
\qquad\qquad\quad\diy\tp (1)\frac{D^{(1)}\gam_1-D^{(2)}\gam_2}{1+D^{(1)}\gam_1-D^{(2)}\gam_2}
-\tp (0)\frac{D^{(1)}\gam_1-D^{(2)}\gam_2}{1-D^{(1)}\gam_1+D^{(2)}\gam_2}=0.\ena\eeq
The last equation (coming from (i3)) can be reduced to two alternative ones:
\beq\label{eq2alt1}{\rm{(A)}}\;\;D^{(1)}\gam_1-D^{(2)}\gam_2=\tp (1)-\tp (0)\;\hbox{ or }\;{\rm{(B)}}\;\;D^{(1)}\gam_1-D^{(2)}\gam_2=0.\eeq
When $1-b\geq\tp (1)-\tp (0)$, these equations specify two (parallel) segments inside the quadrilateral
\beq\label{eqFSet0}\beal
\bbD:=\Big\{\uD=(D^{(1)},D^{(2)}):\;\;D^{(1)},D^{(2)}\geq 0,\\
\qquad\qquad\qquad D^{(1)}+D^{(2)}< 1,\;
\; b-1\leq D^{(1)}\gam_1-D^{(2)}\gam_2\leq 1- b\Big\}.\ena\eeq

Note that $\bet (\uD)$ depends only on $D^{(1)}\gam_1 -D^{(2)}\gam_2$. Further,
function $\uD\in\bbD\mapsto\alph (\uD)$ is concave, and for the gradient vector of $\bet (\uD)$ we have:
$$\beal\nabla\bet (\uD)=\bigg(\diy\frac{\tp(1)\gam_1}{1+D^{(1)}\gam_1 -D^{(2)}\gam_2}
-\frac{\tp(0)\gam_1}{1-D^{(1)}\gam_1 +D^{(2)}\gam_2},\\
\qquad\qquad\qquad\qquad\qquad\qquad\diy -\frac{\tp(1)\gam_2}{1+D^{(1)}\gam_1 -D^{(2)}\gam_2}
+\frac{\tp(0)\gam_2}{1-D^{(1)}\gam_1 +D^{(2)}\gam_2}\bigg),\ena$$
with
$$\nabla\bet (\uD)=\U0\;\hbox{ iff }\;D^{(1)}\gam_1-D^{(2)}\gam_2=\tp (1)-\tp (0),\;
\hbox{ i.e., $\uD$ \ lies in segment (A) in \eqref{eq2alt1}.}$$

We see that when $b\leq 2\tp (0)$, the optimal value, \ $\bet^\rO$, for $\bet (\uD) $ in \eqref{eqOPT1} is attained
at each point of segment (A). Then $n\bet^\rO$ gives the maximal value for $\bbE S_n$
among all strategies $\{\uC_j\}$ satisfying (a0)--(a2) in Eqn  \eqref{eqaVsyo}.

On the other hand, if $2\tp (0)<b<1$ then only the alternative (B) in \eqref{eq2alt1} remains. In this case,
the optimal value in problem \eqref{eqOPT1} is $0$. Thus, $\bbE S_n=0$ yields the maximal value
among all strategies $\{C_j\}$ satisfying (a0)--(a3) in \eqref{eqaVsyo}. This is attained at each point \
$\uD$ \ of segment (B), including $\uD=\U0$.

\section{A general discrete-time setting}

The constructions developed so far show that the Markovian setting is not necessary for the main result.
In this section we attempt to propose a general background where the proposed techniques still works.
Here we again deal with results $\veps_n$
of subsequent random trials, $n=0,1,\ldots$. No specific condition upon the joint distribution is assumed, apart from a dominance condition (involving a given WF); see Eqn
\eqref{eq5.22} below. Each $\veps_n$ is a random element
in a standard measure space $(\cX_n,\fX_n,\mu_n)$. We suppose that a random string
$\beps_0^n=(\veps_0,\ldots ,\veps_n)$ has a joint probability density or probability mass function
(PD/MF) $f_n(\ux_1^n)$ relative to reference measures $\mu_0^n=\prod\limits_{j=0}^n\mu_j$,
on $\operatornamewithlimits{\times}\limits_{j=0}^n\fX_j$:
\beq\label{eq1}
\bbP (\beps_0^n\in A) =\int_Af_n(\ux_1^n)\rd\omu_n (\ux_1^n),\;\ux_1^n=( x_1, \ldots ,x_n)
\in\operatornamewithlimits{\times}\limits_{j=1}^n\cX_j\;\hbox{ and }\;
A\subseteq\operatornamewithlimits{\times}\limits_{j=1}^n\fX_j .
\eeq
A conditional PD/MF $\tf_n(x_n|\ux_1^{n-1})$ will be also used, with
\beq\label{eq2}
\tf_n(x_n|\ux_1^{n-1})f_{n-1}(\ux_1^{n-1})=f_n(\ux_1^n)\;\hbox{ and }\;
\int_{\cX_n}\tf_n(x_n|\ux_1^{n-1})\rd\mu_n(x_n)=1,\,\;f_n\hbox{-a.s.}
\eeq
When necessary, standard properties of measures $\mu_n$ are assumed by default
(completeness, $\sigma$-finiteness).

Next, suppose that a sequence  $\{g_n,\,n\geq 1\}$ of real-valued return functions is given, where
 $g_n: (\ux_1^{n-1}, x_n)\in\cX_1^{n-1}\times\cX_n\mapsto\bbR$. (For $n=1$ we deal with $g_1(x_1)$.)
As above, if you stake $\$\,c_n$ on the $n$th trial, you win $\$\,c_ng_n(\ux_1^{n-1}, x_n)$ if the
outcome is
$x_n\in\cX_n$ preceded by the string $\ux_1^{n-1}$. That is, you make a profit when
$c_ng_n(\ux_1^{n-1}, x_n)>0$
and incur a loss when $c_ng_n(\ux_1^{n-1}, x_n)<0$.$^{*)}$\footnote{$^{*)}
$All functions figuring throughout the paper are assumed measurable, with a specific indication of the
$\sigma$-algebra when necessary.}
%It is supposed that $\int_{\cX_n}{\mathbf 1}(g(x_n)<0)\rd\mu_n(x_n) >0$.

As before, $Z_0>0$ stands for an initial (random) capital.
%Set: $\fW_0=\sigma (Z_0)$ and $\fW_n=\sigma (Z_0,\beps_0^n)$ for $n\geq 1$.
For a given $n\geq 1$, let $Z_n>0$ denote the
capital after $n$ trials. Then the recursion similar to \eqref{eq13B1} emerges:
\beq\label{eq03}
Z_n=Z_{n-1}+C_{n-1}g_n(\beps_0^{n-1}, \veps_n)=Z_{n-1}\left(1+\frac{C_{n-1}g_n
(\beps_0^{n-1}, \veps_n)}{Z_{n-1}}\right),\;\; n\geq 1.\eeq

Formally, we assume that the probability space under consideration is $(\Om ,\fW,\bbP)$
defined as follows. The sample set $\Om $ is the Cartesian product
$\bbR_+\times\left(\operatornamewithlimits{\times}\limits_{n\geq 0}\cX_n\right)$ where $\bbR_+=(0,\infty )$
equipped with the product $\sigma$-algebra
$\fW=\fB(\bbR_+)\times\left(\operatornamewithlimits{\times}\limits_{n\geq 0}\fX_n\right)$ and the filtration $\fW_n=
\fB\times\left(\operatornamewithlimits{\times}\limits_{0\leq j\leq n}\fX_j\right)$. Here $\fB=\fB (\bbR_+)$ is the Borel
$\sigma$-algebra in $\bbR_+$. Adopting the set-up from Sect {\bf II}, $\fW_0=\sigma (Z_0,\veps_0)$,
and $\fW_n=\fW_0\vee\sigma (\beps_0^n)$ for $n\geq 1$.
Then $\om\in\Om$ is a pair $(z_0,\ux_0^\infty )$ where $\ux_0^\infty$ is a sequence
$\{x_n:\,n\geq 1\}$ where $x_n\in\cX_n$. All random elements above and in the sequel are defined as functions
of $\om\in\Om$, subject to standard measurability assumptions.  A probability measure $\bbP$ is given by
a compatible family of joint PD/MF ${\ov f}_n(z_0,\ux_0^n)$ with respect to $\nu\times\mu_0^n$ where
$\nu$ is a chosen measure on $(\bbR_+,\fB)$ (typically, a counting measure on a finite or countable subset or a Lebesgue measure).
Once more we assume that the random variable
$C_{n-1}=C_{n-1}(\beps_0^{n-1})$ is $\fW_{n-1}$-measurable, i.e., yields
a predictable strategy. That is, the stake in the $n$th trial is
based on the results of trials $1$, $\ldots$, $n-1$.   Then $Z_{n-1}=Z_{n-1}
(\beps_0^{n-1})$ is $\fW_{n-1}$-measurable.

As in the previous sections, we consider the weighted logarithmic growth $S_n$ after $n$ trials:
\beq\label{eqlogGR}S_n:=\sum\limits_{j=1}^n\vphi_j(\beps_0^{j-1}, \veps_j)\ln\,\diy\frac{Z_j}{Z_{j-1}}.
\eeq
Cf. Eqn \eqref{eq14B1}. The goal is the same: to maximize, in $\{C_0,\ldots,C_{n-1}\}$, the mean-value
$\bbE S_n$ (and to identify maximizers).
Here we deal with general WFs $\ux_1^j\mapsto \vphi_j(\ux_1^{j-1}, x_j)\geq 0$
depending on the current outcome $x_j$
and the vector of preceding outcomes $\ux_1^{j-1}$.  Again, we can think that
$\vphi_j(\ux_1^{j-1}, x_j)$ represents a `utility' value of outcome $x_j$ (given that it
succeeds an outcome sequence $\ux_1^{j-1}$ from previous trials), and it is taken into
account when we calculate $S_n$.
As above, when $\vphi_j\equiv 1$, the sum \eqref{eqlogGR} becomes telescopic and equal to
$\diy\ln\,\frac{Z_N}{Z_0}$.

The maximization procedure involves a sequence of a.s. positive CFs $\tq_j(\ux_1^{j-1}, x_j)$, $j\geq 0$,
figuring in Eqns \eqref{eq7} -- \eqref{eq5.22} below. Like we said, typically,
the function $q_j$ will be a (conditional) PDF relative to the reference measure $\mu_j$
on $\cX_j$. Define RVs $\alph_j=\alph_j(\beps_0^j)$ and $\rA_n =\rA_n(\beps_0^{n-1})$ by
\beq\label{eq7}
\beacl\diy\alph_j&=\int_{\cX_{j+1}}\vphi_{j+1}(\beps_0^j, x_{j+1})\tf_{j+1}(x_{j+1}|\beps_0^j)
\ln\,\diy\frac{\tf_{j+1}(x_{j+1}|\beps_0^j)}{\tq_{j+1}(\beps_0^j, x_{j+1})}\rd\mu_{j+1}(x_{j+1})\\
\;&\diy =\bbE\;\,\left[\vphi_{j+1}(\beps_0^j, \veps_{j+1})\ln\,\frac{\tf_{j+1}(\veps_{j+1}|
\beps_0^j)}{\tq_{j+1}(\beps_0^j, \veps_{j+1})}\Big|\fW_j\right],\;\hbox{ and \
$\rA_n:=\diy\sum\limits_{j=1}^n\alpha_{j-1}$.}\ena\eeq
Here the RV $\alpha_n$ represents the weighted KL entropy of the conditional
PD/MF $\tf_n(\;\cdot\;|\beps_0^{n-1})$ relative to $q_n (\beps_0^{n-1},\cdot\,)$. Cf. Eqn \eqref{eq17B1}.
The RV  $\rA_n$ yields  the cumulative weighted KL.
% We can treat $\{A_n\}$ is an RP with increments $\alpha_n$ (a weighted KL entropy RP).
%\newpage %\vskip .5 truecm

Fix $b>0$. As before, it is convenient to summarize the assumptions about RVs $\{C_j\}$: \
$\forall$ \ $j\geq 0$,
\beq\label{eqsust3}\beal {\rm{(a0)}}\quad C_j\in\fW_j,\;\;{\rm{(a1)}}\;\; 0\leq C_j< Z_j,\;\;{\rm{(a2)}}\;\;
1+\diy\frac{C_jg_{j+1}(\beps_0^j,\veps_{j+1})}{Z_j}\geq b,\;\hbox{ and}\\
{\rm{(a3)}}\quad C_j(\beps_0^j)\int_{\cX_{j+1}}\vphi_{j+1}(\beps_0^j,x_{j+1})\tq_{j+1}(\beps_0^j,l)
g_{j+1}(\beps_0^j,x_{j+1})\rd\mu_{j+1}(x_{j+1})=0,\ena\eeq
referred to, respectively, as predictability, sustainability, no-ruin and (q,g)-balance conditions.

We also assume that \ $\forall$ \ $j\geq 0$,
\beq\label{eq5.22}\beal
{\rm{(q-f)}}\quad\int_{\cX_j}\vphi_n(\beps_0^{j-1}, x_j )\tq_j (\beps_0^{j-1},
x_j )\rd\mu_j(x_j)\leq\bbE\big[\vphi_j(\beps_0^{j-1},\veps_j)\big|\fW_{j-1}\big] \\
\qquad\qquad\qquad{} =\;\int_{\cX_j}
\vphi_n(\beps_0^{j-1},x_j)\tf_j(x_j|\beps_0^{j-1})\rd \mu_j(x_j)\;\;
\hbox{(weighted $(q,\tf)$-dominance).}\ena\eeq

The conditions (i1)--(i4) are re-written as follows: \\

(i) $\forall$ \ $j\geq 0$, \ $\exists$ a RV $D_j (\beps_0^j)$ with the properties
\beq\label{eq5.0}\beal
{\rm{(i1)}}\;\;\hbox{$0\leq D_j (\beps_0^j)< 1$ (D-sustainability),}\;
{\rm{(i2)}}\;\;\hbox{$1+D_j(\beps_0^j)g_j(\beps_0^j, \veps_{j+1})\geq b$ (D-no-ruin),}\\
{\rm{(i3)}}\;\;D_j(\beps_0^j)\int_{\cX_{j+1}}\diy\frac{\vphi_{j+1}(\beps_0^j, x_{j+1})\tf_{j+1}(x_{j+1}|\beps_0^j)g_{j+1}(\beps_0^j,x_{j+1})}{
1+D_j(\beps_0^j)g_{j+1}(\beps_0^j;x_{j+1})}\rd\mu_{j+1} (x_{j+1})=0\\
\qquad \hbox{(WE (D,g)-balance), with alternatives }\\
{\rm{(i3A)}}\;\int_{\cX_{j+1}}\diy\frac{\vphi_{j+1}(\beps_0^j, x_{j+1})\tf_{j+1}(x_{j+1}|\beps_0^j)g_{j+1}(\beps_0^j, x_{j+1})}{
1+D_j(\beps_0^j)g_{j+1}(\beps_0^j, x_{j+1})}\rd\mu_{j+1} (x_{j+1})=0,{\rm{(i3B)}}\;
D_j(\beps_0^j)=0,\\
\hbox{and }\;\;{\rm{(i4)}}\quad\tq_j(\beps_0^{j-1}, \veps_j)=\diy\frac{\tf_j(\veps_j|\beps_0^{j-1})}{1
+D_{j-1} (\beps_0^{j-1})g_j(\beps_0^{j-1}, \veps_j)}\quad\hbox{ (q-representation)}.
\ena\eeq
As was stressed, Eqns \eqref{eqsust3}--\eqref{eq5.0} are assumed \
$f_j$-a.s. Such a convention is also extended to similar relations
below. All integrals involved are supposed to converge absolutely. As before, we
can interpret these conditions as implemented by a risk-averse trader.

\vskip .5 truecm

{\bf Theorem 4.1.} \label{thm:A} {\sl Assume we are given CFs $\tq_j>0$, WFs
$\vphi_n\geq 0$ and
RFs $g_n$ with values in $\bbR$ satisfying \eqref{eq5.22}. Assume the recursion \eqref{eq03}.
Consider the RVs $S_n$ and $\rA_n$
from Eqns \eqref{eqlogGR} and \eqref{eq7}. Then: \vskip .3 truecm

{\rm{(a)}} For any sequence of RVs $\{C_j,\,j\geq 0\}$ satisfying
\eqref{eqsust3}, the sequence
$\{S_n-\rA_n,\,n\geq 1\}\;$ is a supermartingale. Consequently, $\bbE\,S_n\leq \bbE\,\rA_n$.
\vskip .3 truecm

{\rm{(b)}} To achieve equality $\bbE\,S_n= \bbE\,\rA_n$: the sequence
$\{S_n-\rA_n\}$ is a martingale for a sequence \ $\{C_j\}$ satisfying
\eqref{eqsust3} iff \ $\forall$ \ $j\geq 0$ iff the properties {\rm{(i)}}, {\rm{(ii)}} below
hold true. \medskip

\noindent

\quad{\rm{(i)}} $\forall$ \ $j\geq 0$, \ $\exists$ a RV $D_j (\beps_0^j)$ such that the relations {\rm{(i1)}}--{\rm{(i3)}}
in Eqn \eqref{eq5.0} hold true, and CFs $q_j$ are given by the equation {\rm{(i4)}}.

\quad{\rm{(ii)}} The strategy $C_j$, $j\geq 0$ yields a proportional investment: $C_j (\beps_0^j)=D_j(\beps_0^j)Z_j$.

Furthermore, the CF $q_j$ from {\rm{(i4)}} has $\int_{\cX_j}q_j(\ux_0^{j-1}, x_j)\rd\mu_j(x_j)=1$
(i.e., determines a PD/MF)
iff, in addition to {\rm{(i3)}}, we have that, $\forall$ $j\geq 1$,
\beq\label{eqMBaa}D_{j-1}(\beps_0^{j-1})=0\;\hbox{ or }\;\int_{\cX_j}\frac{\tf_j(x_j|\beps_0^{j-1})
g_j(\beps_0^{j-1},x_j)}{1
+D_{j-1} (\beps_0^{j-1})g_j(\beps_0^{j-1};x_j)}\,\rd\mu_j (x_j)=0.\eeq

{\rm{(c)}} Suppose that \ $\forall$ \ $j\geq 0$, \ $\exists$ \ a RV $D_j(\beps_0^j)$ satisfying
{\rm{(i1)}}--{\rm{(i3)}}. Construct the RVs $D_j^\rO(\beps_0^j)$ in the following manner.
Set $D_j^\rO(\beps_0^j)=D_j(\beps_0^j)$ if the alternative {\rm{(i3A)}} is fulfilled and the value
$D_j(\beps_0^j)>0$ (such a value, if it exists, is unique), and $D_j^\rO(\beps_0^j)=0$ otherwise.
Then the policies $C^\rO_j=D^\rO(\beps_0^j)Z_j$ yield the expectations $\bbE S_n=E_n$ where
\beq\label{eqOptRateG}E_n=\sum\limits_{j=1}^n\bet_{j-1},\;\hbox{with}\; \bet_{j-1}=\bbE\Big\{\vphi_j(\beps_0^{j-1},\veps_j)
\ln\,\Big[1+D^\rO_{j-1}(\beps_0^{j-1})g_j(\beps_0^{j-1},\veps_j)\Big]\Big\}.\eeq
It is a maximal value of $\bbE S_n$ among the strategies satisfying properties {\rm{(a1)}}--{\rm{(a3)}}
in Eqn \eqref{eqsust3}.

{\rm{(d)}} Under the assumptions adopted in {\rm{(c)}}, suppose that \ $\forall$ \ $j\geq 0$, the
RV $D^\rO_j(\beps_0^j)>0$ (so, the alternative {\rm{(i3A)}} is fulfilled with $D_j(\beps_0^j)>0$). Then setting $C^\rO_j=D^\rO_jZ_j$
yields the policies that maximize the mean value $\bbE S_n$ over the strategies satisfying
properties {\rm{(a1)}}--{\rm{(a2)}} in Eqn \eqref{eqsust3}.} \vskip .3 truecm

{\bf Proof of Theorem 4.1}. We still follow the previously established pattern.
(a) Write:
$$\beal\bbE\Big\{\big(S_n-\rA_n\big)\big|\fW_{n-1}\Big\}=S_{n-1}-\rA_{n-1}\diy+\bbE\left\{\left[
\vphi_n\ln\,\left(1+\frac{C_ng_n}{Z_{n-1}}
\right)\right]\Big|\fW_{n-1}\right\}-\alph_n.\ena$$
Next, represent
\beq\label{eq10}
\beal\diy\bbE\left\{\left[\vphi_n (\beps_0^{n-1}, \veps_n)\ln\,\left(1+\frac{C_ng_n(\beps_0^{n-1}, \veps_n)}{Z_{n-1}}\right)\right]\Big|\fW_{n-1}\right\}-\alph_n \\
\quad\diy =\;\int_{\cX_n}\; \vphi_n(\beps_0^{n-1}, x_n)\;\tf_n(x_n|\beps_0^{n-1})\;\ln\;\left[1 +
\frac{C_ng_n(\beps_0^{n-1}, x_n)}{Z_{n-1}}\right] \rd\mu_n(x_n)\\
\qquad\qquad\diy - \int_{\cX_n}\vphi_n(\beps_0^{n-1}, x_n)\tf_n(x_n|\beps_0^{n-1})\ln\,
\frac{\tf_n(x_n|\beps_0^{n-1})}{\tq_n (\beps_0^{n-1}, x_n)}\rd\mu_n(x_n)\\
\quad\diy = \int_{\cX_n}\vphi_n(\beps_0^{n-1}, x_n)\tf_n(x_n|\beps_0^{n-1}) \ln\frac{1
+C_ng_n(\beps_0^{n-1}, x_n)/Z_{n-1}}{\tf_n(x_n|\beps_0^{n-1})/\tq_n (\beps_0^{n-1}, x_n )} \rd\mu_n(x_n)\\
\quad\diy = \int_{\cX_n}\vphi_n(\beps_0^{n-1}, x_n)\;\tf_n(x_n|\beps_0^{n-1})\;
\ln\;\frac{h_n(\beps_0^{n-1}, x_n)}{\tf_n(x_n|\beps_0^{n-1} )}\; \rd\mu_n(x_n)\\
\quad\diy \leq \int_{\cX_n}\vphi_n(\beps_0^{n-1}, x_n)\tf_n(x_n|\beps_0^{n-1})\\
\qquad\qquad\qquad\times\diy\;\left[ \frac{h_n(\beps_0^{n-1}, x_n )}{\tf_n
(x_n|\beps_0^{n-1})} - 1 \right]{\mathbf 1}\big(\tf_n
(x_n|\beps_0^{n-1})>0\big)\; \rd\mu_n(x_n)\\
\quad\diy = \int_{\cX_n}\vphi_n(\beps_0^{n-1}, x_n) \Big[h_n(\beps_0^{n-1}, x_n) - \tf_n(x_n|\beps_0^{n-1})
\Big]\rd\mu_n(x_n) \leq 0.\ena\eeq
Here $h_n(\beps_0^{n-1}, x_n) := \diy\tq_n (\beps_0^{n-1}, x_n )\left[1 + \frac{C_ng_n(\beps_0^{n-1}, x_n )
}{Z_{n-1}}\right]$, $x_n\in\cX_n$. The final inequality in \eqref{eq10} holds since, almost surely,
\beq\label{eq11.1}\beal
\diy \int_{\cX_n}\vphi_n(\beps_0^{n-1}, x_n)h_n(\beps_0^{n-1}, x_n ) \rd\mu_n(x_n) = \int_{\cX_n}
\vphi_n(\beps_0^{n-1}, x_n)\tq_n (\beps_0^{n-1}, x_n ) \rd\mu_n(x_n)\\
\qquad+\diy\frac{C_n}{Z_{n-1}}\int_{\cX_n}\vphi_n(\beps_0^{n-1}, x_n)\tq_n (\beps_0^{n-1}, x_n)
g_n(\beps_0^{n-1}, x_n)\rd\mu_n(x_n)\\
\qquad\qquad\qquad\qquad\qquad\qquad\qquad\qquad\quad\diy \leq \int_{\cX_n}\vphi_n(\beps_0^{n-1}, x_n)
\tf_n(x_n|\beps_0^{n-1})\rd\mu_n(x_n),\ena\eeq
due to (q--f) in \eqref{eq5.22}.

As a result, we get the supermartingale inequality
\beq\label{eq12}
\bbE\left\{\Big[S_n-\rA_n\Big]\Big|\fW_{n-1}\right\}\leq S_{n-1}-\rA_{n-1}.
\eeq

(b) For the martingale property we need to attain equalities in Eqn \eqref{eq10}. The first inequality
becomes equality iff \ $\diy\left[ \frac{h_n(\beps_0^{n-1}, \veps_n )}{\tf_n(\veps_n|\beps_0^{n-1})} - 1 \right]{\mathbf 1}
\big(\tf_n(\veps_n|\beps_0^{n-1})>0\big)=0$ \ $\mu_n$-a.s., \  i.e.,
\beq\label{eq11.2}q_n (\beps_0^{n-1}, \veps_n )\left[1 + \frac{C_ng_n(\beps_0^{n-1}, \veps_n )}{Z_{n-1}}\right]=\tf_n
(\veps_n|\beps_0^{n-1}),\eeq
which yields representation (i4) in Eqn \eqref{eq5.0}. The second equality, achieved in \eqref{eq11.1}, follows from \eqref{eq5.22} after substituting (i4) and is equivalent to (i3). Equation \eqref{eqMBaa} is also
established by using (i4).

Finally, properties (c,d) follows (a) and (b).  $\Box$

\vskip .5 truecm

{\bf Remarks.} {\bf 4.1.} Assumption \eqref{eq5.0} is meaningful when functions
$x_n\in\cX_n\mapsto g_n(\beps_0^{n-1}, x_n)$ are lower-bounded for a.a. $\beps_0^{n-1}$.

{\bf 4.2.} Condition (a1) is not used in the proof of assertion (a), only (a2) is relevant there.

{\bf 4.3.} Taking $q_n(x_n, \ux_1^{n-1})=
\tf_n(x_n|\ux_1^{n-1})$ leads to the case $S_n=0$.

{\bf  4.4.} As earlier, the staple of the proof of Theorem 4.1 is the Gibbs
inequality for weighted entropies; see \cite{SSYK, SYS}. It is similar to the standard Gibbs
inequality (cf. \cite{CT, KS}) but requires additional assumptions \eqref{eq5.22} and \eqref{eq5.0}.

{\bf 4.5.} The optimization problem emerging in assertion (c) is an example of
a control problem. The connections with specific facts and methods, including the Bellman
equation seem fruitful and should be explored in forthcoming works.

\vskip .5 truecm

The next step is to provide an assertion covering Scheme II (a riskless asset
and a collection of risky assets) in a general setting. The trial RP $\{\veps_n,n\geq 0\}$
has the same (general) nature as before, and we continue using the same concepts and
notations whenever possible. By the time of the $n$th trial, we now have $K(n)+1$ assets
under management,
which generates  RFs $g^{(0)}_n$ and $\ug_n=\left(g^{(1)}_n, \ldots, g^{(K(n))}_n\right)$  at
the $n$th trial. We assume that the RF $g^{(0)}$ is (relatively) riskless and has the form
$$\beal g^{(0)}_n(\ux_0^{(n-1)})=1+\rho_{n-1}\;\hbox{ where }\;\rho_{n-1}
=\rho_{n-1} (\ux_0^{(n-1)})\geq 0\\
\qquad\qquad\hbox{(an interest rate between the times of the $(n-1)$st and the $n$th trial).}\ena$$
The remaining RFs are risky, with values $g^{(s)}_n(\ux_0^{(n-1)},x_n)$, $1\leq s\leq K(n)$,
depending on the outcomes $\ux_0^{(n-1)}$ of the previous trials (which you know) and
the outcome $x_n$ of the $n$th trial (which is unknown at the time of the $n$th investment).
We again set $\u1_n=(1,\ldots ,1)$ (with $K(n)$ entries altogether) and
\beq\label{eq414}\ug^*_n=\ug_n-(1+\rho_{n-1})\u1_n=\left(g^{(1)}_n-(1+\rho_n),\ldots ,g^{(K(n))}_n-(1+\rho_n)\right).\eeq

Similarly to \eqref{eqRL&sR}, the recursion for the capital value is taken to be
\beq\label{eqRL&sRG}\beacl Z_n&=(1+\rho_{n-1})(Z_{n-1}-|C_{n-1}|)+\uC_{n-1}\cdot
\ug^*_n(\veps_{n-1},\veps_n)\\
\;&\qquad\qquad\qquad\qquad\qquad\diy =Z_{n-1}\left[
1+\rho_{n-1} +\diy\frac{\uC_{n-1}\cdot\ug^*_n(\beps_0^{n-1},\veps_n)}{Z_{n-1}}\right],\;n\geq 1.\ena\eeq
Here the vector $\uC_j=\left(C^{(1)}_j,\ldots ,C^{(K(n))}_j\right)$ represents a
random portfolio which determines the investment at the time after the $j$th trial, and
$\uC_j\cdot\ug^*_{j+1}(\veps_j,\veps_{j+1})=\sum\limits_{s=1}^{K(j)}C_{j-1}^{(s)}
\left[g_{j+1}^{(s)}(\beps_0^{j-1},\veps_j)-(1+\rho_j)\right]$.
The norm $|\uC_j|=\sum\limits_{s=1}^{K(j)} C^{(s)}_{j}=\uC_j\cdot\u1_j$ shows the total size of
the investment.
Given a sequence of CFs $(\ux_0^{j-1},x_j)\mapsto q_j(\ux_0^{j-1},x_j)>0$ and the threshold
value $b>0$, we work under the following assumptions: \ $\forall$ \ $j\geq 0$,
\beq\label{eqa127}\beal{\rm{(a0)}}\;\;\uC_j\in\fW_j,\;\;{\rm{(a1)}}\;\;\uC_j\geq\U0,\;|\uC_j|< Z_j,\\
{\rm{(a2)}}\;1+\rho_j(\beps_0^j)+\diy\frac{\uC_j(\beps_0^j)\cdot\ug^*_{j+1}(\beps_0^j,\veps_{j+1})}{
Z_j(\beps_0^j )}\geq b,\;\hbox{ and}\\
{\rm{(a3)}}\;\; \int_{\cX_{j+1}}\vphi_{j+1}(\veps_j,x_{j+1})\tq_{j+1}(\veps_0^j,x_{j+1})
\left[\uC_j(\beps_0^j)\cdot
\ug^*(\veps_j,x_{j+1})\right]\rd\mu_{j+1}(x_{j+1})=0.\ena\eeq

As before, we want to maximize the mean value $\bbE S_n$. Here, as before,
$$S_n:=\sum\limits_{j=1}^n
\vphi_j(\beps_0^{j-1},\veps_j)\ln\diy\frac{Z_j}{Z_{j-1}}$$
and $(\ux_0^{j-1},x_j)\in \cX_0^{j-1}\times \cX_j\mapsto\vphi_j(\ux_0^{j-1},x_j)\geq 0$ is
a given sequence of WFs $\vphi_j$, $j\geq 1$. The maximization problem is  in
$\{\uC_j,\,j\geq 0\}$ subject to restrictions (a0)--(a2) or (a0)--(a3) in \eqref{eqa127}.

As usual, we assume the (q,f)-dominance property, that $\forall$ $j\geq 0$,
\beq\label{eqq12RLG}\beac{\rm{(q-f)}}\;\;
\int_{\cX_{j+1}}\vphi_{j+1} (\beps_0^j,x_{j+1})\Big\{\left[1+\rho_j(\beps_0^j)\right]
\q_{j+1}(\beps_0^j,x_{j+1})-\tf_{j+1} (x_{j+1}|\beps_0^j)\Big\}\leq 0.\ena\eeq

Next, consider modified conditions (i1)--(i4):

\quad \ (i) $\forall$ \ $j\geq 0$ \ $\exists$ a random vector  $\uD_j(\beps_0^j)$ where
$\uD_j(\beps_0^j)=\left(D^{(1)}_j(\beps_0^j),\ldots ,D^{(K(n))}_j(\beps_0^j)\right)$, such that
\beq\label{eq8M137}\beal\quad  {\rm{(i1)}}\; \uD_j(\beps_0^j)\geq\U0,\;\hbox{ and }\;
|\uD_j(\beps_0^j)|:=\sum\limits_{s=1}^{K(j)}D^{(s)}_j(\beps_0^j)<1,\\
\quad {\rm{(i2)}}\; 1+\rho_j(\beps_0^j)+\left[\uD_j(\beps_0^j)\cdot\ug^*_{j+1}(\beps_0^j,\veps_{j+1})\right]
\geq b,\\
\quad {\rm{(i3)}}\;\int_{\cX_{j+1}}\tf_{j+1}(x_{j+1}|\beps_0^j)\,
\diy\frac{\vphi_{j+1}(\beps_0^j,x_n)\left[\uD_j(\beps_0^j)\cdot\ug^*_{j+1}(\beps_0^j,x_{j+1})\right]}{1+\rho_j(\beps_0^j)
+\uD_j(\beps_0^j)\cdot\ug^*_{j+1}(\beps_0^j,x_{j+1})}\rd\mu_{j+1}(x_{j+1})= 0,\\
\quad {\rm{(i4)}}\;\hbox{ the CFs have the form }\;\tq_{j+1}(\beps_0^j,\veps_{j+1})=\diy\frac{\tf_{j+1}(\veps_{j+1}|\beps_0^j)}{1
+\rho_j(\beps_0^j) +\uD_j (\beps_0^j)\cdot\ug^*_{j+1}(\beps_0^j,\veps_{j+1})}.\ena\eeq

The behavior of $\bbE S_n$ is characterized through the cumulative weighted KL entropy RP
$\{\rA_n,\,n\geq 1\}$ where $\rA_n$ and $\alph_j (\veps_{j-1})$ are  as in Eqn \eqref{eq7}.

\vskip .5 truecm

{\bf Theorem 4.2.} {\sl Adopt the above setting in Eqns \eqref{eq414}--\eqref{eq8M137}. Then
the following assertions hold true.

{\rm{(a)}} Take any sequence of RVs $\uC_j=
\left(C^{(1)}_n,\ldots ,C^{(K(j))}_j\right)\in\fW_j$, $j\geq 0$, satisfying \eqref{eqa127}. Then
the sequence $\{S_n-\rA_n,\,n\geq 1\}\;$ is a supermartingale; hence $\bbE\,S_n\leq \,
\bbE \rA_n$ \ $\forall$ \ $n\geq 1$.

{\rm{(b)}} For equality $\bbE\,S_n=\sum\limits_{j=1}^n\bbE \alph (\veps_{j-1})$: the sequence
$\{S_n-\rA_n,\,n\geq 1\}$ is a martingale for predictable RVs \ $\uC_j$, $j\geq 0$, satisfying
\eqref{eqa127} \ iff \ the additional properties {\rm{(i)}}, {\rm{(ii)}} below are fulfilled. \vskip .3 truecm

\quad{\rm{(i)}} $\exists$ \ a sequence of RVs \ $\uD_j(\beps_0^j)=\left(D^{(1)}_j(\beps_0^j),\ldots ,
D^{(K(j))}_j(\beps_0^j)\right)$ \ such that \ $\forall$ \ $j\geq 0$, properties
{\rm{(i1)}}--{\rm{(i4)}} in Eqn \eqref{eq8M137} hold true.

\quad{\rm{(ii)}} $\forall$ \ $j\geq 0$, \ the policy $\uC_j$ produces a proportional investment:
after the $j$th trial
the amount $\uC_j(\beps_0^j)=\uD_j(\beps_0^j)Z_j$ goes to the assets with
the RF $\ug_j$ whereas $Z_j-|D_{j}|Z_{j}$ is allocated to the riskless return.
\medskip

{\rm{(c)}} Suppose that \ $\forall$ \ $j\geq 0$, \ $\exists$ a RV $\uD_j(\beps_0^j)=
\left(D^{(1)}_j(\beps_0^j),\ldots ,D^{(K(j))}_j(\beps_0^j)\right)$
fulfilling conditions {\rm{(i1)}}--{\rm{(i3)}} in Eqn \eqref{eq8M137}. Let $\ubD$ stand
for a sequence of RVs $\{\uD_j\}$. Then, for
the corresponding proportional investment portfolios $\uC_j$, the expectation
$E_n(\ubD)=\bbE S_n$ has the form
\beq\label{eqOptRateGr}\beal E_n(\ubD )=\sum\limits_{j=0}^{n-1}\bet_{j-1}\;\hbox{ where }\;
\bet_{j-1}=\bbE\Big\{\vphi (\veps_{j-1},\veps_j)
\ln\,\Big[1+\rho +\uD_{j-1}(\beps_0^{j-1})\cdot\ug^*(\beps_0^{j-1},\veps_j)\Big]\Big\}.\ena\eeq
Then the maximum of $E_n(\ubD )$ over the RVs satisfying {\rm{(i1)}}--{\rm{(i3)}}
yields the maximum of the mean value $\bbE S_n$ over the strategies satisfying
{\rm{(a0)}}--{\rm{(a3)}}  in Eqn \eqref{eqa127}.

{\rm{(d)}} Suppose that \ $\forall$ \ $j\geq 0$, \ $\exists$ a RV $\uD_j(\beps_0^j)=
\left(D^{(1)}_j(\beps_0^j),\ldots ,D^{(K(j))}_j(\beps_0^j)\right)$ with all entries
$D^{(s)}_j(\beps_0^j)>0$, fulfilling conditions {\rm{(i1)}}--{\rm{(i3)}} in Eqn \eqref{eq8M137}.
Then the maximum of $E_n(\ubD )$ from Eqn \eqref{eqOptRateGr} yields the maximum
the mean value $\bbE S_n$ among the strategies satisfying conditions
{\rm{(a0)}}--{\rm{(a2)}}  in Eqn \eqref{eqa127}. }

\medskip

We omit the proof of Theorem 4.2 as it does not contain new elements compared with
the preceding statements.

\vskip .5 truecm

{\bf Example 4.1: Uniform IID trials.} As an illustration, consider again the case of IID trials,
within the remits of Scheme I, with a single asset.
Suppose $\cX_n=\bbR$ and $\mu_n$ is the Lebesgue measure. Take
$\tf_n (x_n|\ux_0^{n-1})=f(x_n)$, i.e., assume the results of the trials are independent and distributed
with PDF $f$: $\veps_n\sim f$.  Take a WF $\vphi (x)\geq 0$, a bounded RF $g(x)$
and a CF $q(x)>0$. Also fix $b\in (0,1)$. As in the above IID examples, the
weighted KL entropy $\alph_j$ from Eqn \eqref{eq7} becomes a constant:
$$\beac\alph_{j-1}=\alph =\int\vphi (x)f(x)\ln\diy\frac{f(x)}{q(x)}\rd x=\bbE\left[\vphi (\veps_1)
\ln\frac{f(\veps_1)}{q(\veps_1)}\right].\ena$$
The martingale condition (i4) in Theorem 4.1(b)
is that $q(x)=\diy\frac{f(x)}{1+Dg(x)}$ where, in accordance with (i1) and (i2) from
Eqn \eqref{eq5.0}, $D\in [0,1]$ is a constant such that $1+Dg(x)\geq b$,
for $f$-a.a. \ $x\in\bbR$. In addition, (i3) requires that
$$ D\int {\diy\frac{f(x)\vphi (x)g(x)}{1+Dg(x)}}\rd x=D\,\bbE\diy\frac{ \vphi (\veps_1)g(\veps_1)}{1
+Dg(\veps_1)}= 0.$$

We see that, following \eqref{eqOptRateG}, the function
$$\bet (D)=\int\vphi (x)f(x)\ln [1+Dg(x)]\rd x=\bbE\vphi (\veps_1)\ln [1+Dg(\veps_1)]$$
should be maximised in the variable $D$, subject to the above restrictions (i1)--(i3).
Set:
$$D_+=\max\;\Big\{D\in [0,1]:\;1+Dg(x)\geq b\;\forall\;x\in{\rm{supp}}\;f\Big\}.$$
Then
$D\in [0,D_+]\mapsto \bet (D)$ is a concave function, with the derivative ${\diy\frac{\rd}{\rd D}}\bet (D)
=\bbE{\diy\frac{\vphi (\veps_1)g(\veps_1)}{1+Dg(\veps_1)}}$.
Thus, we are interested in a stationary point $D_0$, where $\bbE{\diy\frac{\vphi (\veps_1)g(\veps_1)}{1+D_0g(\veps_1)}}=0$. At $D=0$ we have $\bet (D)=0$ and
${\diy\frac{\rd}{\rd D}}\bet (D)=\int\vphi (x)f(x)g(x)\rd x:=\bbE\vphi (\veps_1)g(\veps_1)$. We see that if
$\bbE\vphi (\veps_1)g(\veps_1)= 0$ then the optimal proportion $D^\rO$ from Theorem 4.1(c)
equals $0$, i.e., a cautious trader would not invest
in the IID market where on average there is no profit/loss. (Let alone the case where the loss
exceeds the profit.) Otherwise, i.e.,
when  $\bbE\vphi (\veps_1)g(\veps_1)>0$, the trader looks at the value $D_0$: if $D_0\leq D_+$,
we set $D^\rO=D_0$. Otherwise, when $D_0> D_+$, the theory does not give a formal answer.

\vskip .5 truecm

To be more specific, consider an IID case with uniformly distributed IID trial results,
where \\
$f\sim$U$[a_1,a_2]$, with $a_1<0<a_2$. Thus, $f(x)=\B1(a_1\leq x\leq a_2)\big/a$
where $a=a_2-a_1$. Take $\vphi (x)\equiv 1$ (no preference); examples of a non-constant
WF can be also incorporated into the argument that follows. Fix constants $b>0$ and
$\delta_\pm ,\gamma_\pm\in\bbR$ and consider a piece-wise linear RF
$$g(x)=\bcs \delta_+x+\gam_+,&0< x\leq a_2,\\ \delta_-x +\gam_-,&a_1\leq x< 0.
\ecs$$
For definiteness, assume that $g(x)\geq 0$ for $0\leq x\leq a_2$ and
$g(x)\leq 0$ for $a_1\leq x\leq 0$; this boils down to $\gam_+\geq 0$ and $a_2\delta_+\,
+\gam_+\geq 0$ and the opposite inequalities for $\delta_-$ and $\gam_-$ with $a_1$ replacing $a_2$.
(The case $\delta_\pm  =0$ was effectively treated in Example 1.1.)  For the mean of the RF we
obtain
$$\bbE g(\eps )=\diy\frac{1}{a}\left[a_2\gam_++a_1\gam_-+\frac{a_2^2\delta_+}{2}
-\frac{a_1^2\delta_-}{2}\right].$$

To guarantee the martingale condition,  we have to take $q(x)$ as in (i4):
$q(x)=\diy\frac{f(x)}{1+Dg(x)}$ where $D$ obeys (i1)--(i2): $0\leq D\leq D_+$
where
$$D_+=\max\,\Big[D\in [0,1]:\,1+D\gam_-\geq b\;\hbox{ and }
1+D\big(a_1\delta_-+\gam_-\big)\geq b\Big].$$
Then, referring to \eqref{eqOptRateG}, we want to maximise in $D$ the function
$D\mapsto\bet (D)$. Here
$$\bet (D)={\diy\frac{1}{a}}\int_{a_1}^0\ln\big\{1+D[x\delta_-+\gam_-]\big\}\rd x
+{\diy\frac{1}{a}}\int_0^{a_2}\ln\big\{1+D[x\delta_++\gam_+]\big\}\rd x,$$
with $\bet (0)=0$. A direct calculation yields
$$\frac{\rd}{\rd D}\bet (D)=\frac{1}{a}\left(\int_{a_1}^0\frac{(\gam_-+x\delta_-)\rd x}{
1+D\gam_- +xD\delta_-}
+\int_0^{a_2}\frac{(\gam_++x\delta_+)\rd x}{1+D\gam_++xD\delta_+}\right)$$
with $\diy\frac{\rd}{\rd D}\bet (D)\big|_{D=0}=\bbE g(\eps )$.
This allows us to specify the above routine in terms of $a_{1,2}$, $\delta_\pm$ and $\gam_\pm$.
Namely,
if $\bbE g(\eps )\leq 0$ then the optimal policy is $D^{\rO}=0$ (no investment),
while if $\bbE g(\eps )>0$ then $D^\rO=D_0$ whenever $D_0\leq D_+$.
Here $D_0>0$ is the zero of the derivative
$\diy\frac{\rd}{\rd D}\bet (D)$ solving a transcendental equation
$$a-\frac{1}{D_0\delta_+}\ln\frac{1+D_0\gam_+}{1+D_0\gam_++D_0a_2\delta_+}
+\diy\frac{1}{D_0\delta_{-}}\ln\frac{1+D_0\gam_-}{1+D_0\gam_-+D_0a_1\delta_-}=0.$$
As was stated earlier, for $D_0> D_+$ the theory does not yield a formal answer.

A similar methodology is applicable in examples where the PDF $f$ and RF $g$ are
piece-wise polynomial functions.

\vskip .5 truecm

{\bf Example 4.2: Gaussian IID trials.} Another example is where the trial results
$\veps_n$ are IID and have a normal PDF: $f\sim$N$(0,\sigma^2)$. Here
\beq\label{eq19}
f(x)=\frac{\exp\,(-x^2/\sigma^2)}{\sigma {\sqrt{2\pi\,}}},
\;x\in\bbR.\eeq
Adopt Scheme I, again with a single asset. We have a choice
of RFs $g$ and WFs $\vphi$ where one can do explicit calculations. Viz, take
$$g(x)=\bcs a_1x+a_2,&x>0,\\ -a_3,&x<0,\ecs$$
where $a_1,a_2\geq 0$, $a_3> 0$.
Then the WF can be taken piece-wise polynomial, e.g.,
$$\vphi (x)=\bcs x^2\theta_++x\gam_++\delta_+,&x>0,\\
x^2\theta_-+x\gam_-+\delta_-,&x<0.\ecs$$
The mean value $\bbE \vphi (\veps_1)g(\veps)$ can be computed as a polynomial
in variables $a_{1,2,3}$, $\theta_\pm$, $\gam_\pm$ and $\delta_\pm$.

Let us also fix $b\in (0,1)$. For the martingale CF $q(x)=\diy\frac{f(x)}{1+Dg(x)}$, the weighted
KL entropy $\bet =\bet (D)$ becomes
$$\bet =\ln\,(1-Da_3)\int_{-\infty}^0f(x)\vphi (x)\rd x+
\int_0^\infty f(x)\vphi (x)\ln\,[1+D(a_1x+a_2)]\rd x,$$
which is defined for $0\leq D\leq D_+:=1\wedge [(1-b)/a_3]$. Here the derivative
$$\frac{\rd}{\rd D}\bet (D)=-\frac{a_3}{1-Da_3}\int_{-\infty}^0f(x)\vphi (x)\rd x+
\int_0^\infty\frac{ f(x)\vphi (x)(a_1x+a_2)}{1+D(a_1x+a_2)}\rd x,$$
again with $\diy\frac{\rd}{\rd D}\bet (D)\Big|_{D=0}=\bbE \vphi (\veps_1)g(\veps_1)$.
The argument from the previous example can be still used to specify the optimal
proportion $D^\rO$.

\vskip .5 truecm

{\bf Example 4.3. IID trials with linear/logarithmic risky RFs.} In this example, we assume
the trial results $\veps_n$ are IID and uniformly distributed on the interval $[-1,1]$. Thus,
$\cX_n=[-1,1]$, the reference
measures
$\mu_n$ are Lebesgue's, and the PDF $\tf_n(x_n|\ux_0^{n-1})=\diy\frac{1}{2}$ \ $\forall$ $\ux_1^{n-1}
\in [0,1]^{n-1}$ and $x_n\in [0,1]$. Adopt Scheme II and
suppose that we have one riskless asset with the RF $g^{(0)}(x)
=1+\rho$ and two risky assets, with
$$\beac\hbox{$g^{(1)}(x)=-\gam x$ and $g^{(2)}(x)=
-\theta\ln \,(1-x)$,}\\
\ug^*(x)=(-\gam x-(1+\rho),-\theta\ln\,(1-x)-(1+\rho)),\ena\;-1\leq x\leq 1.$$
Here $\rho \geq 0$,  $\gam >0$ and $\theta >0$ are given parameters.
(Additional bounds for $\rho$, $\gam$ and $\theta$ will appear as simplifying assumptions below.)
Also fix $b>0$ and a WF $x\in [-1,1]\mapsto\vphi (x)\geq 0$.

The inequalities (i1), (i2) in \eqref{eq8M137} take now the form
$$\beac{\rm{(i1)}}\;\;D^{(1)},D^{(2)}\geq 0\;\hbox{ and }\;D^{(1)}+D^{(2)}<1 \;\hbox{(D-sustainability)},\\
{\rm{(i2)}}\;\;(1+\rho)(1-D^{(1)} -D^{(2)})-D^{(1)}\gam x-D^{(2)}\theta\ln\,(1-x)\geq b,\;-1\leq x\leq 1\;\hbox{(D-no-ruin).}\ena$$
The function $x\in [-1,1]\mapsto (1+\rho)(1-D^{(1)}-D^{(2)}) -D^{(1)}\gam x -D^{(2)}\theta\ln\,(1-x)$ in the LHS of (i2)
has the second derivative $\diy\frac{D^{(2)}\theta}{(1-x)^2}\geq 0$, so is convex. Its minimum
on $[-1,1]$ is attained at $x=1-\diy\frac{D^{(2)}\theta}{D^{(1)}\gam}$ \ if \
$D^{(2)}\theta\leq 2D^{(1)}\gam$ and at $x=-1$ if $D^{(2)}\theta >2D^{(1)}\gam$.
So, we can re-write the condition in the form
\beq\label{eqFinEx1}{\rm{(i2)}}\;\;\bcs (1+\rho)(1-D^{(1)}-D^{(2)})-D^{(1)}\gam \left(1-\diy\frac{D^{(2)}\theta}{D^{(1)}\gam}\right)
-D^{(2)}\theta\ln\,\diy\frac{D^{(2)}\theta}{D^{(1)}\gam}\geq b,&\hbox{if }\;
D^{(2)}\theta\leq 2D^{(1)}\gam ,\\
(1+\rho)(1-D^{(1)}-D^{(2)})+D^{(1)}\gam-D^{(2)}\theta\ln\, 2\geq b,&\hbox{if }\;D^{(2)}\theta \geq 2D^{(1)}\gam .\ecs\eeq

For definiteness, let us focus on the bottom line in \eqref{eqFinEx1} and consider the following
domain $\bbD$ in the
$D^{(1)},D^{(2)}$-plane:
$$\beal\bbD =\Bigg\{\uD =\left(D^{(1)},D^{(2)}\right):\;D^{(1)},D^{(2)}\geq 0,\;D^{(1)}+D^{(2)}< 1,\\
\qquad\qquad\qquad\qquad\qquad 2D^{(1)}\gam\leq D^{(2)}\theta\leq\diy\frac{(1+\rho)(1-D^{(1)}) -b+D^{(1)}\gam}
{\ln\,2+(1+\rho)/\theta}\Bigg\}.\ena$$
Suppose that
$$0\vee 1+\rho-\ln\,2-\frac{1+\rho}{\theta}<b<1+\rho.$$
Then $\bbD$ is non-empty and represents a triangle or  quadrilateral with vertices at
$\uD=\U0 =(0,0)$ (the origin) and
$\uD =\left(0,\diy\frac{1+\rho -b}{\ln\,2+(1+\rho)/\theta}\right)$ (a point in the $D^{(2)}$ -axis);
the remaining two or one lie(s) on  the straight lines where $2D^{(1)}\gam = D^{(2)}\theta$
and $D^{(2)}\theta =\diy\frac{(1+\rho)(1-D^{(1)}) -b+D^{(1)}\gam}{\ln\,2+(1+\rho)/\theta}$.

Next, consider condition (i3) in \eqref{eq8M137}:
\beq\label{eqCondi3FE}
\int_{-1}^1\frac{\vphi (x)\Big\{D^{(1)}(\gam x+1+\rho )+D^{(2)}[\theta\ln\,(1-x)+1+\rho]\Big\}}{
1+\rho -D^{(1)}(\gam x+1+\rho )-D^{(2)}[\theta\ln\,(1-x)+1+\rho]}\rd x=0\;
\hbox{(WE D,g-balance)}.\eeq
In the strong form it becomes a system of two equations, for the maximum
of the concave function $(D^{(1)},D^{(2)})\mapsto\bbE\diy\ln\,[1+\rho+\uD\cdot\ug^*(\veps_1)]$:
\beq\label{eqStCndi3FE}\beal\diy\int_{-1}^1\frac{\vphi (x)(\gam x+1+\rho )}{1+\rho -
D^{(1)}(\gam x+1+\rho )-D^{(2)}[\theta\ln\,(1-x)+1+\rho]}\rd x\\
\qquad\diy =
\int_{-1}^1\frac{\vphi (x)[\theta\ln\,(1-x)+1+\rho]}{1+\rho -
D^{(1)}(\gam x+1+\rho )-D^{(2)}[\theta\ln\,(1-x)+1+\rho]}\rd x=0.\ena\eeq

Obviously, $\uD=\U0$, i.e., $D^{(1)}=D^{(2)}=0$ is a solution to \eqref{eqCondi3FE}.
If it is the only solution in $\bbD$ then we have $\uD^\rO=\U0$ (no investment in the risky assets, with
$\bbE S_n=n\ln\,(1+\rho )$). It yields
the optimum over the portfolio sequences $\uC_j=(C^{(1)}_j,C^{(2)}_j)$ such that, \ $\forall$ $j\geq 0$,
\beq\label{eqa197}\beal{\rm{(a0)}}\;\;\uC_j\in\fW_j,\;\;{\rm{(a1)}}\;\;C^{(1)}_j(\beps_0^j),C^{(2)}_j(\beps_0^j)\geq 0,
\;\;C^{(1)}_j(\beps_0^j)+C^{(2)}_j(\beps_0^j)< Z_j,\\
{\rm{(a2)}}\;\;1+\rho -\diy\frac{C^{(1)}_j(\beps_0^j)(\gam x+1+\rho )
+C^{(2)}_j(\beps_0^j)[\theta\ln\,(1-x)+1+\rho ]}{
Z_j(\beps_0^j)}\geq b,\;-1\leq x\leq 1,\\
\hbox{and }\;{\rm{(a3)}}\;\;\diy\int_{-1}^1\frac{\vphi (y)\big\{C^{(1)}_j(\beps_0^j)(-\gam y-1-\rho)+
C^{(1)}_j(\beps_0^j)[-\theta\ln\,(1-y)-1-\rho]\big\}}{1+\rho+
C^{(1)}_j(\beps_0^j)(-\gam y-1-\rho)+
C^{(1)}_j(\beps_0^j)[-\theta\ln\,(1-y)-1-\rho]} \rd y=0.\ena\eeq

An alternative possibility may occur if we have a solution to \eqref{eqStCndi3FE}
lying in $\bbD$. In this case we set $\uD^\rO$ to be the solution of \eqref{eqStCndi3FE}
and achieve the maximum of $\beta (\uD )=\bbE\ln\,[1+\rho \uD\cdot\ug^*(\veps_1)]$ in $\bbD$ at $\uD=\uD^\rO$.
Then $\bbE S_n=n\beta (\uD^\rO)$ yields the optimum among all portfolio sequences
$\uC_j=(C^{(1)}_j,C^{(2)}_j)$ satisfying, \ $\forall$ $j\geq 0$, the requirements (a0)--(a2), without
(a3). Note that still $\uD^\rO=\U0$ if we have that
$$\diy\int_{-1}^1\vphi (x)(\gam x+1+\rho )\rd x =
\int_{-1}^1\vphi (x)[\theta\ln\,(1-x)+1+\rho]\rd x=0.$$

%\newpage
\vskip .5 truecm

\noindent
{\bf Acknowledgement.}
The authors are grateful to their respective institutions for the financial support and
hospitality: MK thanks the Higher School of Economics for
the support in the framework of a subsidy granted to the HSE by the Government of the Russian Federation for the implementation of the Global Competitiveness Programme; IS was supported by the
FAPESP Grant No 11/51845-5, and thanks the IMS, USP and to Math Dept, University of
Denver; YS thanks the Math Dept, PSU.

\end{document}